\documentclass[11pt,reqno]{amsart}

\usepackage[top=100pt,bottom=60pt,left=96pt,right=92pt]{geometry}

\usepackage{hyperref}
\usepackage{amsmath,amsthm,amssymb,amsfonts,enumerate,xcolor,layout,tikz,subcaption,enumitem,fullpage, verbatim, colordvi}
\usetikzlibrary{calc}   
\usepackage{stmaryrd}
\usepackage{graphicx}
\usepackage{tikz}
\usepackage{pgfplots}
\usepackage{subcaption}
\pgfplotsset{compat=1.18}
\usetikzlibrary{arrows.meta}

\newtheorem{maintheorem}{Theorem}
\newtheorem{theorem}{Theorem}[section]
\newtheorem*{theorem*}{Theorem}
\newtheorem{proposition}[theorem]{Proposition}
\newtheorem{corollary}[theorem]{Corollary}

\newtheorem{lemma}[theorem]{Lemma}

\theoremstyle{definition}
\newtheorem{definition}[theorem]{Definition}
\newtheorem{remark}[theorem]{Remark}
\newtheorem{example}[theorem]{Example}

\theoremstyle{plain}

\definecolor{linkblue}{rgb}{0,0,.6}
\definecolor{citered}{rgb}{.7,0,0}
\hypersetup{colorlinks =true, linkcolor=linkblue, citecolor = citered, urlcolor=linkblue}

\def\R{{\mathbb R}}
\def\T{{\mathbb T}}

\newcommand{\interior}{\operatorname{Int}}

\newcommand{\purple}[1]{\textbf{\purple{#1}}}

\begin{document}
\title{Fiber-connectivity of almost-toric systems in four-dimension }

\author{Aerim Si}
\address{Department of Mathematics\\
	University of Illinois at Urbana-Champaign\\Champaign, Illinois 61801}
\curraddr{}
\email{aerimsi2@illinois.edu}
\thanks{}

\date{\today}

\begin{abstract}
    In this paper, we study the fiber connectivity of proper almost-toric systems $\Phi$ on connected four-dimensional symplectic manifolds. Using Morse–Kirwan theory, we establish three sufficient conditions for fiber connectivity. These results generalize and extend the fiber connectivity results of Pelayo, Ratiu, and V\~u Ng\d{o}c \cite{pelayo2015fiber} by relaxing some of their hypotheses.
\end{abstract}

\maketitle

\tableofcontents

\section{Introduction}\label{sec:intro}

A fundamental question in  the study of integrable systems is whether the fibers of a moment map are connected. Atiyah \cite{atiyah1982convexity} and Guillemin–Sternberg \cite{GuSt1982} independently proved that the moment map of a Hamiltonian torus action on a compact connected symplectic manifold has connected fibers. This result was extended to compact connected symplectic orbifolds by Lerman–Tolman \cite{lerman1997hamiltonian}, and to connected symplectic orbifolds with proper moment maps by Lerman–Meinrenken–Tolman–Woodward \cite{lerman1996non}. For semitoric systems in dimension four, V\~u Ng\d{o}c  \cite{VN2007} proved that the fibers of a semitoric system are connected; this was  generalized to dimension six by Wacheux \cite{wacheux2013systemes}. More recently, Sepe and Tolman \cite{sepe2026connectedness} generalized fiber-connectedness results for toric and semitoric systems to  complexity one integrable systems with nondegenerate singularities. In subsequent work \cite{sepe2025connectedness}, they further extended these results to allow certain degenerate singularities.

\vspace{1mm}

For almost-toric systems, however, fibers need not be connected in general (see, for instance, \cite[Example 4.10]{pelayo2015fiber}). Nevertheless,  Pelayo, Ratiu, and V\~u Ng\d{o}c \cite[Theorems 1–3]{pelayo2015fiber} established fiber connectivity under additional hypotheses on the moment map image. One of the key ingredients in their approach is Morse–Bott theory. They compose the moment map $\Phi$ with a Morse function $f$ defined on a neighborhood of $\Phi(M)$ and give conditions under which $f\circ \Phi$ is a proper Morse-Bott function with neither index 1 nor coindex 1 (see Theorem 5.6, Proposition 5.8 of \cite{pelayo2015fiber}). Proposition 4.5 of \cite{pelayo2015fiber} then implies that the level sets of $f\circ \Phi$ are connected. 

\vspace{1mm}

In this paper, we generalize and extend known fiber connectivity results for proper almost-toric systems by relaxing some of the hypotheses in \cite[Theorems 1--3]{pelayo2015fiber}.
Our strategy replaces the Morse-Bott framework of \cite{pelayo2015fiber}
 with Morse-Kirwan theory \cite{holm2016morse, kirwan1985cohomology}.
 Let $\Phi:M\rightarrow \mathbb{R}^2$ be a proper almost-toric system on a connected four dimensional symplectic manifold $(M,\omega)$. Fix $c\in \Phi(M)$, we consider the norm-squared function $f_c(x,y)=\left\|(x,y)-c\right\|^2$ and study the composition $h_c:=f_c\circ \Phi=\|\Phi-c\|^2$. We give sufficient conditions for $h_c$ to be a proper, bounded below Morse-Kirwan function  with no critical points of index 1. Since $(h_c)^{-1}(0)=\Phi^{-1}(c)$,  Lemma \ref{main lemma}, which is a Morse-Kirwan extension of  Lemma 4.3 of \cite{pelayo2015fiber}, implies  that $\Phi^{-1}(c)$ is connected. Since  $c\in \Phi(M)$ was arbitrary,  all fibers of $\Phi$ are connected.

\vspace{1mm}

In contrast to \cite[Proposition 4.5]{pelayo2015fiber}, our approach only requires excluding critical points of index 1, but not coindex 1. Using Morse--Kirwan theory, we allow \emph{degenerate} contact (see \cite[Definition 5.5]{pelayo2015fiber}) between the level sets of $f_c$ and the elliptic-regular critical value branches of $\Phi$, provided that every critical point of the restriction of $f_c$ to each branch is a \emph{local minimum}. We also allow critical values of $\Phi$ to be critical points of $f_c$, provided that they are \emph{local minima} of $f_c$. Moreover, by choosing $f_c$ to be proper, we remove the assumption in Theorem 2 \cite{pelayo2015fiber} that $\Phi(M)$ lies in a proper convex cone.

\medskip
\noindent\textbf{Notation.}
Throughout this paper, $(M,\omega)$ is a connected four-dimensional symplectic manifold. We denote by $\text{Crit}(\Phi)$ and $\text{Crit}_1(\Phi)$ the sets of critical points and rank 1 critical points of $\Phi$, respectively.

\medskip

We can now state our main results. 

\begin{maintheorem}\label{theorem 1}
    Let $\Phi:M\to\mathbb{R}^2$ be a proper almost-toric system on $(M,\omega)$. Suppose there
exists a diffeomorphism $g\colon \Phi(M)\rightarrow\mathbb{R}^2$
onto its image such that $g(\Phi(M))$ is closed in $\mathbb{R}^2$ and, for
every connected component $C\subset \operatorname{Crit}_1(\Phi)$, $g\bigl(\Phi(C)\bigr)$
is contained in an affine line. 
Then all fibers of $\Phi$ are connected.
\end{maintheorem}
In particular, Theorem \ref{theorem 1} covers new examples beyond the criteria of \cite[Theorem 1--3]{pelayo2015fiber}, including systems whose moment map images are closed infinite strips, closed half-planes, or the entire $\mathbb{R}^2$ with only focus--focus singularities.

\begin{maintheorem}\label{theorem 2}
    Let $\Phi\colon M \to \mathbb{R}^2$ be a proper almost-toric system on $(M,\omega)$. Suppose there exists a diffeomorphism $g\colon \Phi(M)\rightarrow\mathbb{R}^2$ onto its image such that, for $\widetilde{\Phi}:=g\circ \Phi$, the image $\widetilde{\Phi}(M)$ is a closed convex subset of $\mathbb{R}^2$ and every elliptic-regular critical value of $\widetilde{\Phi}$  lies in $\partial \widetilde{\Phi}(M)$. Then all fibers of $\Phi$ are connected. 
\end{maintheorem}
Theorem \ref{theorem 2} not only recovers the fiber connectivity results for all four cases considered in \cite[Theorem 3]{pelayo2015fiber} but also extends them to noncompact symplectic manifolds.

\begin{maintheorem}\label{theorem 3}
    Let $\Phi\colon M \to \mathbb{R}^2$ be a proper almost-toric system on $(M,\omega)$. Suppose there exists a diffeomorphism $g\colon \Phi(M)\rightarrow\mathbb{R}^2$ onto its image such that, for $\widetilde{\Phi}:=g\circ \Phi$, the image $\widetilde{\Phi}(M)$ is closed in $\mathbb{R}^2$,  every elliptic-regular critical value of $\widetilde{\Phi}$  lies in $\partial\widetilde{\Phi}(M)$, and 
    \[\partial \widetilde{\Phi}(M)=\bigsqcup_{i=1}^n \Gamma_i, \quad \Gamma_i=R_i(\Gamma), \quad n\geq 3,\]
    where $R_i$ is the rotation about the origin through the angle $\theta_i=\frac{2\pi(i-1)}{n}$ for $i=1,\cdots,n$, and 
    \[\Gamma=\{(x,y) \; | \; \frac{x^2}{\alpha^2}-\frac{y^2}{\beta^2}=1, \ x>0    \}, \quad  \alpha, \beta >0, \quad \frac{\beta}{\alpha}=tan(\frac{\pi}{n}).\] Then all fibers of $\Phi$ are connected. 
\end{maintheorem}


\newcommand{\hyperbolicplot}[1]{%
\begin{tikzpicture}
\begin{axis}[
    width=\linewidth,
    height=\linewidth,
    xmin=-5, xmax=5,
    ymin=-5, ymax=5,
    axis equal image,
    axis lines=none,
    ticks=none,
    enlargelimits=false,
    clip=true,
    samples=180
]

\pgfmathsetmacro{\B}{tan(180/#1)}
\pgfmathsetmacro{\T}{3.5}


\addplot[
    draw=none,
    fill=gray!20,
    domain=0:360,
    samples=900,
    variable=\a
]
(
 {cos(\a)/
  sqrt(
    cos(mod(\a+180/#1,360/#1)-180/#1)^2
    -sin(mod(\a+180/#1,360/#1)-180/#1)^2/
     (\B*\B)
  )},
 {sin(\a)/
  sqrt(
    cos(mod(\a+180/#1,360/#1)-180/#1)^2
    -sin(mod(\a+180/#1,360/#1)-180/#1)^2/
     (\B*\B)
  )}
) -- cycle;

\foreach \k in {0,...,\numexpr#1-1\relax} {
    \pgfmathsetmacro{\ang}{360*\k/#1}

    \foreach \sgn in {-1,1} {
        \pgfmathsetmacro{\rayang}
            {\ang+\sgn*180/#1}

        \addplot[
            gray!60,
            dashed,
            thin,
            domain=0:8,
            variable=\r
        ]
        ({\r*cos(\rayang)},
         {\r*sin(\rayang)});
    }
}

\foreach \k in {0,...,\numexpr#1-1\relax} {
    \pgfmathsetmacro{\ang}{360*\k/#1}

    \addplot[
        black!75,
        thick,
        domain=-\T:\T,
        variable=\t
    ]
    (
      {cosh(\t)*cos(\ang)
       -\B*sinh(\t)*sin(\ang)},
      {cosh(\t)*sin(\ang)
       +\B*sinh(\t)*cos(\ang)}
    );
}

\end{axis}
\end{tikzpicture}%
}

\begin{figure}[htbp]
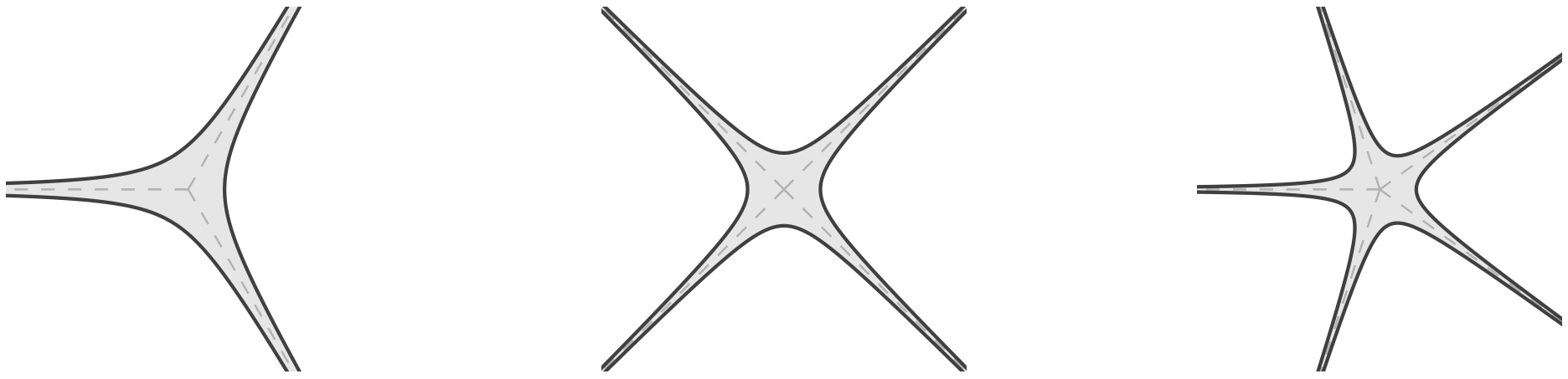

\centering

\begin{subfigure}[b]{0.27\textwidth}
    \centering
    \hyperbolicplot{3}
    \caption*{$n=3$}
\end{subfigure}%
\hspace{0.015\textwidth}%
\begin{subfigure}[b]{0.27\textwidth}
    \centering
    \hyperbolicplot{4}
    \caption*{$n=4$}
\end{subfigure}%
\hspace{0.015\textwidth}%
\begin{subfigure}[b]{0.27\textwidth}
    \centering
    \hyperbolicplot{5}
    \caption*{$n=5$}
\end{subfigure}

\caption{
Image of $\widetilde{\Phi}$ described in Theorem \ref{theorem 3} for $n=3,4,5$, with $\alpha=1$ and
$\beta=\tan(\pi/n)$.
The solid curves represent the rotated
hyperbolic branches $\Gamma_i$,
and the dashed rays indicate their asymptotes.
}
\label{fig:hyperbolic-boundaries}
\end{figure}

Theorem \ref{theorem 3} provides a further class of almost-toric systems with connected fibers whose moment map images are noncompact and nonconvex, with boundaries consisting of disjoint branches of hyperbolas. Together, these results provide new geometric criteria for fiber connectivity beyond the hypotheses considered in \cite{pelayo2015fiber}. We expect that our approach will also be useful for studying fiber connectivity of higher-dimensional almost-toric systems.

\medskip
\noindent\textbf{Structure of the paper.}
In Section \ref{sec:background}, we review background on nondegenerate integrable systems, their local normal forms, and Morse-Kirwan functions. In Sections \ref{sec:ER}--\ref{sec:FF}, we analyze the composition $f\circ \Phi$ near elliptic-regular, elliptic-elliptic, and focus-focus singular points, respectively. In Section \ref{sec:6}, we combine these local results to give conditions under which $f\circ \Phi$ is a Morse-Kirwan function with no points of index 1. Finally, in Section \ref{sec:main}, we apply these  results to norm-squared functions on $\mathbb{R}^2$ to establish criteria for fiber connectivity and prove our main theorems.

\medskip
\noindent\textbf{Acknowledgements.} The author is deeply grateful to her advisor, Professor Susan Tolman, for suggesting the use of Morse–Kirwan theory to extend previous fiber-connectivity results, and for her invaluable guidance and support throughout this project. The author also sincerely thanks Professors Yael Karshon and Joseph Palmer, and Yichen (Jason) Liu, for insightful discussions that contributed to the development of this paper.

\section{Background}\label{sec:background}

In  this section, we briefly recall the basic notions of integrable systems and their nondegenerate singularities, together with the local normal form. For further background, we refer  to \cite{bolsinov2004integrable, leung2010almost, palmer2025integrable, pelayo2011symplectic, sepe2018integrable}. We also review Morse-Kirwan functions and the results needed later; see \cite{holm2016morse, kirwan1985cohomology} for further details.

\medskip
\noindent\textbf{Integrable systems, nondegnerate singularities, and local normal form.}   
\begin{definition}
    A \emph{completely integrable system} on  $(M,\omega)$ consists of smooth functions $\Phi_1,\Phi_2 \in C^{\infty}(M)$ such that  $\{\Phi_1,\Phi_2\}=0$, and whose differentials $d \Phi_1,d \Phi_2$ are linearly independent almost everywhere. These functions define a map $\Phi:=(\Phi_1,\Phi_2):M\xrightarrow{}\mathbb{R}^2$,  which is often called the moment map.
\end{definition}

\begin{definition}
An integrable system $\Phi=(\Phi_1,\Phi_2):M\xrightarrow{} \mathbb{R}^2$ is called \emph{almost-toric} if all of its singularities are non-degenerate with no hyperbolic components.  
\end{definition}

To describe the singularities that occur in almost-toric systems, we  recall the notion of nondegenerate critical point of integrable systems.

\begin{definition}
    Let $\Phi=(\Phi_1,\Phi_2)$ be an integrable system on $(M,\omega)$ and let $p\in M$ be a critical point of $\Phi$. If $\text{rank}(d\Phi)_p=0$, $p$ is called \emph{nondegenerate} if the Hessians of $\Phi_1,\Phi_2$ span a Cartan subalgebra of the Lie algebra of quadratic forms on $T_pM$.  If $\operatorname{rank}(d\Phi)_p=1$, after a linear change of the components, assume that $d\Phi_1(p)\neq 0,\; d\Phi_2(p)=0.$ Let $L=\mathbb{R}X_{\Phi_1}(p)\subset T_pM$ be the one-dimensional
subspace spanned by the Hamiltonian vector field $X_{\Phi_1}$ at $p$, and let $L^\omega$ denote the symplectic orthogonal complement of $L$. Then $p$ is called
\emph{nondegenerate} if $\operatorname{Hess}_p(\Phi_2)$ on $L^\omega/L$ is a
nondegenerate symmetric bilinear form.
\end{definition}

The space $Q(T_pM)$ of quadratic forms on $T_pM$, endowed with the Poisson bracket induced by $\omega_p$, is a Lie algebra naturally isomorphic to $\mathfrak{sp}(T_pM,\omega_p)$. Williamson \cite{williamson1936algebraic} classified its Cartan subalgebras in terms of elliptic, hyperbolic, and focus-focus components. This gives rise to the Williamson classification of nondegenerate rank zero critical points (see Theorem 1.3 and 1.6 of \cite{bolsinov2004integrable}). The Williamson classification is reflected in the local normal forms of nondegenerate critical points. The following local normal forms is due to Eliasson \cite{Eli84, eliasson1990normal}, V\~u Ng\d{o}c-Wacheux \cite{vu2013smooth}, and Miranda-Zung\cite{miranda2006note}. 

\begin{theorem}\label{local normal form theorem}(Local normal form)
    Let $p\in M$ be a nondegenerate critical point of an integrable system $\big(M,\omega, \Phi=(\Phi_1,\Phi_2) \big)$. Then  there exist local symplectic coordinates $(x,\xi)=(x_1,x_2,\xi_1,\xi_2)$ on an open neighborhood $U\subset M$ of $p$ such that $x_i(p)=\xi_i(p)=0$, and a map $q=(q_1,q_2):U\rightarrow \mathbb{R}$ such that $\{\Phi_i,q_j \}=0$ for all $i,j=1,2$ whose components $q_i$ are of the following forms: 
    \begin{itemize}
            \item Regular component: $q_i=\xi_i$,
            \item Elliptic component: $q_i=\frac{x_i^2+\xi_i^2}{2}$,
            \item Focus-focus component: $q_1=x_1\xi_2-x_2\xi_1$,  $q_2=x_1\xi_1+x_2\xi_2$,
            \item Hyperbolic component: $q_i=x_i\xi_i$.
        \end{itemize}
        Moreover, if $p$ has no hyperbolic component, then there exists a local diffeomorphism $g:\mathbb{R}^2 \rightarrow \mathbb{R}^2$ with $g(0)=\Phi(p)$ such that $\Phi(x,\xi)=(g\circ q)(x,\xi)$ in a neighborhood of $p$.

\end{theorem}

Thus, for an almost-toric system in dimension four, a rank zero critical point is either of elliptic-elliptic or focus-focus type, while a rank one critical point is of elliptic-regular type. 

We next give two connectivity lemmas  that will be used later. The first gives a general criterion for fiber connectivity, while the second concerns the local connectivity of regular values of an almost toric system.

\begin{lemma}\label{detour lemma}(\cite{pelayo2015fiber}, Lemma 4.6) Let $f:M\to\mathbb{R}^n$ be a smooth map from a smooth connected
manifold $M$ to $\mathbb{R}^n$. Let $B_r$ be the set of regular values
of $f$. Suppose that $f$ has the following properties.
\begin{enumerate}
    \item $f$ is a proper map.
    
    \item For every sufficiently small neighborhood $D$ of any critical
    value of $f$, $B_r\cap D$ is connected.
    
    \item The regular fibers of $f$ are connected.
    
    \item The set $\operatorname{Crit}(f)$ of critical points of $f$
    has empty interior.
\end{enumerate}
Then the fibers of $f$ are connected.
\end{lemma}

\begin{lemma}\label{locally connectedness of B_r}
    Let $\Phi\colon M\to\mathbb{R}^2$ be a proper almost-toric integrable system on $(M,\omega)$
and let $B:=\Phi(M)$. Suppose that $B$ is convex and  the set $B_r$ of regular values of $\Phi$ satisfies $B_r=\interior{(B)}\setminus\mathcal{F}$ where $\mathcal{F}$ is the set of focus-focus critical values of $\Phi$. Then,  for any $c\in B$ and for any sufficiently small disc $D$ centered at $c$, the set $B_r\cap D$ is connected.
\end{lemma}

\begin{proof}
    Let $D$ be an open disc centered at $c \in B$. Since  both $D$ and $\interior{(B)}$ are open convex subsets of  $\mathbb{R}^2$, their intersection $U:=D\cap \interior{(B)}$ is also an open convex subset of $\mathbb{R}^2$. By assumption, $B_r\cap D=\big(\interior{(B)}\setminus\mathcal{F} \big)\cap D=U \setminus (\mathcal{F}\cap D)$. Since $\Phi$ is proper and focus-focus critical points are isolated, the set $\mathcal{F}$ is locally finite. Thus, since  $\overline{D}$ is compact, $\mathcal{F}\cap \overline{D}$, and hence $\mathcal{F}\cap D$, is finite. Therefore, $B_r \cap D$ is an open convex subset of $\mathbb{R}^2$ with finitely many points removed. Thus $B_r \cap D$ is path connected, so it is connected. \end{proof}

\medskip

\noindent\textbf{Morse--Kirwan functions}
We now recall the notion of a Morse-Kirwan function and some properties, including a key lemma, that will be used throughout the paper. 

\begin{definition}\label{def:Morse-Kirwan}
Let $M$ be a smooth manifold, and let $H:M\to\mathbb{R}$ be a smooth function. Let $C\subset \operatorname{Crit}(H)$ be a closed subset on which $H$ is constant. A submanifold $N\subset M$
containing $C$ is a \emph{minimizing manifold} for $H$ along $C$ if: 
\begin{enumerate}
    \item The tangent space $T_xN$ at every point $x\in C$ is maximal among subspaces of $T_xM$ on which  $\operatorname{Hess}_x H$ is positive semidefinite.
    \item The restriction $ H|_N$ of $H$ to $N$ achieves its minimum exactly on $C$.
\end{enumerate}

The function $H$ is \emph{Morse--Kirwan} or \emph{minimally degenerate} if its critical set is a locally finite disjoint union of closed subsets $\operatorname{Crit}(H)
=\bigsqcup_{\alpha\in A} C_\alpha$ and there exists a minimizing manifold $N_\alpha$ for $H$ along each $C_\alpha$.
\end{definition}

\begin{remark}
    Every Morse--Bott function is Morse--Kirwan.
\end{remark}

\begin{remark}
    The critical subset $C_\alpha$ in Definition \ref{def:Morse-Kirwan} need not be connected. Moreover, a minimizing manifold $N_\alpha$ may have connected components of different dimensions (see \cite[Section 10]{kirwan1985cohomology}),  \cite[Section 2]{holm2016morse}). However, by replacing each $C_\alpha$ with its intersections with the connected components of $N_\alpha$, we may assume that the index of $\operatorname{Hess}_p (H)$ is constant along any $p \in C_\alpha$. This index is equal to the codimension of  $N_\alpha$, and is called the index of $H$ along $C_\alpha$.
\end{remark}

\begin{remark}
Kirwan's original definition \cite{kirwan1985cohomology}  of minimally degenerate  function on a compact manifold requires that the critical set is a \emph{finite} disjoint union of closed subsets and that the corresponding minimizing manifolds are co-orientable. Following Holm--Karshon \cite{holm2016morse}, we adopt the more general definition in which the critical set is  required only to admit  a \emph{locally finite}  disjoint union of closed subsets, and no co-orientability assumption is imposed on the minimizing manifolds. Since $M$ is not assumed to be compact, this general definition is the appropriate one for our setting. On a compact manifold, local finiteness is equivalent to finiteness so two definitions agree
(see Remarks 2.5 and 2.6 of \cite{holm2016morse}). Throughout this paper, a \emph{Morse–Kirwan} function means a minimally degenerate function in the sense of  Definition~\ref{def:Morse-Kirwan}.
\end{remark}

\begin{example}
    Consider \(f:\mathbb R^2\to\mathbb R\) defined by \(f(x,y)=x^4\). Then $C:=\operatorname{Crit}(f)=\{(0,y)\mid y\in\mathbb R\}$. The function \(f\) is not Morse–Bott since \(\operatorname{Hess}_p(f)\) is degenerate in the normal direction to \(C\) for every \(p\in C\), while it is Morse–Kirwan since \(N=\mathbb R^2\) is a minimizing manifold for \(f\) along \(C\). This example shows that Morse–Kirwan functions allow degenerate local minima.
\end{example}

\medskip

The following theorem of Holm--Karshon shows that the Morse--Kirwan condition is local. Consequently, to prove that a function is globally Morse--Kirwan, it suffices to verify the Morse--Kirwan condition in a neighborhood of each point.

\begin{theorem}[\cite{holm2016morse}, Theorem 2.8]\label{local to global}
    Let $M$  be a smooth manifold and $f:M\rightarrow \mathbb{R}$ a smooth function. Suppose that every point in $M$ has an open neighborhood $U$ such that $f|_{U}$ is minimally degenerate. Then $f$ is mininally degenerate.
\end{theorem}

We next give three preliminary results that will be useful for verifying the Morse–Kirwan condition.

\begin{proposition}\label{equivalence of locally finitness}
Let $H \colon M \to \mathbb{R}$ be smooth and let $x \in \operatorname{Crit}(H)$ with $H(x)=0$. The following are equivalent:
\begin{enumerate}
    \item There exists a neighborhood $U$ of $x$ such that $\operatorname{Crit}(H)\cap U$ is a locally finite disjoint union of closed subsets of $U$ on each of which
    $H$ is constant.

    \item There exists a neighborhood $V$ of $x$ such that $\operatorname{Crit}(H)\cap V
        \subset H^{-1}(0).$
\end{enumerate}
\end{proposition}
\begin{proof}
$(1)\Rightarrow(2)$: Assume there exists a neighborhood $U \ni x$ such that
$\text{Crit}(H)\cap U \;=\; \bigsqcup_{\alpha} C_\alpha$, 
where $\{C_\alpha\}$ is a locally finite family of closed subsets of $U$ on each of which $H$ is constant.
Let $A := \{\alpha \ | \  H|_{C_\alpha} \neq 0\}.$ By local finiteness, there exists a neighborhood $W\subset U$ of $x$ such that $F :=\{\alpha \in A \ | \ C_\alpha\cap W \neq\emptyset \}$ is finite. Hence $V := W \setminus \bigcup_{\alpha \in F} C_\alpha$ is an open neighborhood of $x$. Since $V$ does not intersect any $C_\alpha$ with 
$\alpha \in A$, we have
$\operatorname{Crit}(H)\cap V \subset H^{-1}(0).$

$(2)\Rightarrow(1)$: Assume there exists a neighborhood $V$ of $x$ such that $\operatorname{Crit}(H)\cap V \subset H^{-1}(0).$ Take $U := V$ and $C := \text{Crit}(H)\cap V$. Since $\text{Crit}(H)$ is closed in $M$, the set $C$ is closed in $U$. Moreover, $H \equiv 0$ on $C$, since $C\subset H^{-1}(0)$. Also, since $\text{Crit}(H)\cap U=C$, the critical set in $U$ is already a decomposition into a single closed subset on which $H$ is constant. Since a family consisting of a single set is automatically locally finite, it follows that $\text{Crit}(H)\cap U$ is a locally finite disjoint union of closed subsets on which $H$ is constant. 
\end{proof}

\begin{lemma}\label{local min and preimage}
    Let $H\colon M \to \mathbb{R}$ be smooth and let $x \in M$ with $H(x)=0$. Suppose $x$ is a local minimum of $H$ and that there exists a neighborhood $V$ of $x$ such that $\operatorname{Crit}(H)\cap V
        \subset H^{-1}(0).$ Then $H$ is Morse--Kirwan in a neighborhood of $x$.
\end{lemma}
\begin{proof}
    After shrinking $V$ if necessary, we may assume that $H\geq0$ on $V$. Let $C:=\text{Crit}(H)\cap V$. By assumption, $C\subset H^{-1}(0)\cap V$. Since $H\geq0$ on $V$, every point of $H^{-1}(0)\cap V$ is a local minimum of $H$, and hence is critical. Therefore $C=H^{-1}(0)\cap V$. 

    Now note that $H|_{V}$ attains its minimum exactly on $C$. Also, since every $p\in C$ is a local minimum of $H$, $\text{Hess}_p(H)$ is positive semidefinite on $T_pM=T_pV$. Hence $V$ is a minimizing manifold for $H$ along $C$. By Proposition $\ref{equivalence of locally finitness}$, $\text{Crit}(H)\cap V$ satisfies the locally finite decomposition condition in Definition \ref{def:Morse-Kirwan}.  Consequently, $H$ is Morse-Kirwan on $V$.
\end{proof}

\begin{proposition}\label{only local min}
    Let $f:\R \rightarrow \R$ be a smooth function such that every critical point of $f$ is a local minimum. Then $\operatorname{Crit}(f)$ has at most one connected component, and $f$ is a Morse-Kirwan function.  
\end{proposition}
\begin{proof}
    If $\text{Crit}(f)=\emptyset$, then the conclusion is immediate, so suppose $\text{Crit}(f)\neq \emptyset$. We claim that if $a,b$ are critical points of $f$ with $a<b$, then $[a,b] \subset \text{Crit}(f)$.
    
    Let \[ M:=\max_{x\in[a,b]}f(x) \quad \text{and} \quad S:=\{x\in [a,b] \ | \ f(x)=M\}.\] Then $S$ is nonempty and closed in $[a,b]$. If $x\in S\cap (a,b)$, then $x$ is a local maximum of $f$, hence a critical point. Therefore $x$ is also a local minimum by hypothesis. This implies that $f$ is constant near $x$, so $(x-\epsilon,x+\epsilon)\subset S$ for some small enough $\epsilon>0$. 
    If $a\in S$, then $f\leq f(a)$ on $ [a,b]$, while local minimum property of $a$ gives $f\geq f(a)$ near  $a$. Hence $f=f(a)$ on  $[a,a+\epsilon)$ for some $\epsilon>0$, so $[a,a+\epsilon)\subset S$. Similarly, if $b\in S$, then  $(b-\epsilon,b]\subset S $ for some $\epsilon>0$.
    Therefore, $S$ is both open and closed in $[a,b]$, so $S=[a,b]$. Consequently, $f$ is constant on $[a,b]$, so $[a,b]\subset \text{Crit}(f)$. This implies that whenever $p,q\in \text{Crit}(f)$, the entire interval between $p$ and $q$ also lies in $\text{Crit}(f)$. Thus $\text{Crit}(f)$ is connected, and hence has at most one connected component.

    Now, let $C:=\text{Crit}(f)$, and take $N:=\R$. Since every point of $C$ is a local minimum, $\text{Hess}_x(f)$ is positive semidefinite on $T_x\R$ for every $x\in C$.  Hence $T_xN$ is maximal among subspaces on which the Hessian is positive semidefinite. Also, since $C$ is connected and consists entirely of local minima, $f|_N=f$ attains its global minimum exactly on $C$.  Thus $N$ is a minimizing manifold for $f$ along $C$. Since $C$ is closed and locally finite, $f$ is therefore Morse--Kirwan.  \end{proof}
\begin{remark}\label{only local min remark}
    The same proof applies when $\mathbb{R}$ is replaced by an interval, possibly with boundary.
\end{remark}

We next give two conditions under which the composition of a smooth map with a Morse–Kirwan function is Morse–Kirwan.

\begin{lemma}\label{either a regular or an isolated point}
Let $\Phi:X\to Y$ be a smooth map between smooth manifolds, and let
$f:Y\to\mathbb{R}$ be a Morse-Kirwan function. Suppose that  $p\in X$ is either a regular point or an \emph{isolated} critical point of $\Phi$, and that $\Phi(p)$ is a local minimum of $f$. Then $f\circ\Phi$ is a Morse-Kirwan function near $p$.
\end{lemma}

\begin{proof}
    Let $H := f \circ \Phi$. Without loss of generality, we can assume that $H(p) = f(\Phi(p)) = 0$. Since $\Phi(p)$ is a local minimum of $f$, there exists a neighborhood $V$ of $\Phi(p)$ such that $f\geq0$ on $V$. After shrinking a neighborhood $U$ of $p$ so that $\Phi(U)\subset V$, we have $H\geq0$ on $U$. Hence, $p$ is a local minimum of $H$.

We first claim that $\text{Crit}(H)\cap U=\Phi^{-1}(\text{Crit}(f))\cap U$.
If \(p\) is a regular point of \(\Phi\), then after shrinking \(U\) if necessary, \(d\Phi_x\) has rank 2 for all \(x\in U\). Hence, $dH_x = df_{\Phi(x)}\circ d\Phi_x = 0$ if and only if $df_{\Phi(x)}=0$, so the claim follows. If \(p\) is an isolated critical point of \(\Phi\), then after shrinking \(U\) if necessary, we have \(\text{Crit}(\Phi)\cap U=\{p\}\). Thus  \(d\Phi_x\) has rank 2 for all \(x\in U\setminus\{p\}\), and the same argument shows
$\text{Crit}(H)\cap (U\setminus\{p\})=\Phi^{-1}(\text{Crit}(f))\cap (U\setminus\{p\})$. At $p$, since $\Phi(p)$ is a local minimum of $f$, we have \(df_{\Phi(p)}=0\), and therefore \(p\in \text{Crit}(H)\) and \(\Phi(p)\in \text{Crit}(f)\). Thus the equality also holds at $p.$

Since $f$ is Morse--Kirwan, Proposition \ref{equivalence of locally finitness} gives, after shrinking $V$ if necessary, $\text{Crit}(f)\cap V \subset f^{-1}(0)$. Shrinking $U$ so that $\Phi(U)\subset V$, we have 
\[\operatorname{Crit}(H) \cap U =\Phi^{-1}(\operatorname{Crit}(f) ) \cap U 
= \Phi^{-1}(\operatorname{Crit}(f) \cap V) \cap U 
\subset \Phi^{-1}(f^{-1}(0)) = H^{-1}(0).\] Hence, by Lemma \ref{local min and preimage}, $H$ is Morse--Kirwan near $p$.  

\end{proof}

\begin{proposition}\label{submersion and MK}
     Let $\Phi:X\rightarrow Y$ be a submersion, and let $f:Y \rightarrow \mathbb{R}$ be a Morse-Kirwan function. Then $f\circ\Phi:X\rightarrow\mathbb{R}$ is a Morse-Kirwan function. 
\end{proposition}

\begin{proof}
Let $H:=f\circ \Phi$. Since $d\Phi_x$ is surjective  for every $x \in X$, $dH_x = df_{\Phi(x)} \circ d\Phi_x=0$ if and only if $df_{\Phi(x)}=0$. Thus, $\text{Crit}(H)=\Phi^{-1}(\text{Crit}(f))$. Since $f$ is a Morse-Kirwan function, its critical set $\text{Crit}(f)$ admits a locally finite disjoint decomposition $\text{Crit}(f)=\bigsqcup_{\alpha} C_{\alpha}$ where each $C_\alpha$  is closed and $f$ is constant on $C_\alpha$, and there exists a  minimizing manifold $N_\alpha$ for $f$ along $C_\alpha$ for each $\alpha$.

For each $\alpha$, define $\widetilde{C}_\alpha := \Phi^{-1}(C_\alpha)$ and $\widetilde{N}_\alpha := \Phi^{-1}(N_\alpha)$. Since $\Phi$ is continuous and local finiteness is preserved under pullback, the family $\{\widetilde{C}_\alpha\}$ is locally finite and each $\widetilde{C}_\alpha$ is closed.  Thus, $\text{Crit}(H)=\bigsqcup_{\alpha} \widetilde{C}_{\alpha}$ forms a locally finite disjoint union of closed subsets of $H$. Moreover, since $f$ is constant on each $C_{\alpha}$, the function $H$ is constant on each $\widetilde{C}_{\alpha}$. Also, since $N_\alpha\subset Y$ is a submanifold and $\Phi$ is a submersion, $\widetilde{N}_\alpha$ is a submanifold of $X$ and $T_x\widetilde{N}_\alpha=d\Phi_x^{-1}(T_{\Phi(x)}N_\alpha)$.

We now verify that $\widetilde{N}_\alpha$ is a minimizing manifold for $H$  along $\widetilde{C}_\alpha$.  Let $x\in \widetilde{C}_\alpha$ and $y=\Phi(x)$. Since $df_y=0$, we have {\small\begin{align}\label{eqn:hessian}
    \operatorname{Hess}_x(H)(v,w)
=
\operatorname{Hess}_{y}(f)
\bigl(d\Phi_x(v),d\Phi_x(w)\bigr)
+
\underbrace{df_{y}}_{=0}
\bigl(\operatorname{Hess}_x\Phi(v,w)\bigr)
=\operatorname{Hess}_{y}(f)
\bigl(d\Phi_x(v),d\Phi_x(w)\bigr)
\end{align}}for all $v,w\in T_xX$.
Since $T_{y}N_\alpha$ is maximal among subspaces of $T_yY$ on which $\text{Hess}_y(f)$ is positive semidefinite, the Hessian identity (\ref{eqn:hessian}), together with the surjectivity of $d\Phi_x$, implies that 
$d\Phi_x^{-1}(T_yN_\alpha)=T_x\widetilde{N}_\alpha$ is maximal among subspaces of $T_xX$ on which $\text{Hess}_x(H)$ is positive semidefinite. Also, since $f|_{N_\alpha}$ attains its minimum exactly on $C_\alpha$, it follows that
$H|_{\widetilde{N}_\alpha}=f\circ \Phi|_{\widetilde{N}_\alpha}$ attains its minimum exactly on $\widetilde{C}_\alpha$. Therefore, $\widetilde{N}_\alpha$ is a minimizing manifold for $H$ along $\widetilde{C}_\alpha$. Consequently, $H$ is Morse-Kirwan. 
\end{proof}

The following key lemma is the Morse–Kirwan extension of \cite[Lemma 4.3]{pelayo2015fiber}. Its proof follows by generalizing their argument to the Morse-Kirwan setting, using Kirwan's Morse stratification for Morse-Kirwan functions (see \cite[Section 10]{kirwan1985cohomology}) in place of the Morse–Bott stable manifolds.

\begin{lemma}\label{main lemma}
    Let $M$ be a connected manifold and let $f:M\rightarrow \mathbb{R}$ be a proper, bounded from below Morse-Kirwan function with no points of index 1. Then the critical set of index 0 is connected.
\end{lemma}

\begin{proof}
 Let $C_0:=\{p\in \text{Crit}(f)\;|\; \text{index}_p(f)=0\}$. Write $\operatorname{Crit}(f)=\bigsqcup_{\alpha} C_\alpha$
for the locally finite decomposition into critical subsets, and let $\Sigma_\alpha$ be a minimizing manifold for $f$ along $C_\alpha$. Recall that $C_\alpha$ need not be connected and that $\Sigma_\alpha$ may have connected components of different dimensions; write $\Sigma_\alpha=\bigsqcup_i\Sigma_{\alpha,i}$. Thus, the index need not be constant on $C_\alpha$. Following Kirwan (see  \cite[Section 10]{kirwan1985cohomology}),  refine each critical subset $C_\alpha$ by replacing it with intersections $C_{\alpha,i}:=C_\alpha\cap \Sigma_{\alpha,i}$.  Since the original family $\{C_\alpha\}$ is locally finite and the connected components of  $\Sigma_\alpha$ are open in $\Sigma_\alpha$, the refined family $\{C_{\alpha,i}\}$ is also locally finite.
Moreover $\Sigma_{\alpha,i}$ is a minimizing manifold for $f$ along $C_{\alpha,i}$, and  $\text{index}_p(f)=\text{codim}(\Sigma_{\alpha,i})$ for all $p\in C_{\alpha,i}$. Hence the index is constant on each refined critical subsets. After relabeling, we may  assume that the index is constant on each $C_\alpha$, and denote it by  $\text{index}(C_\alpha)$.

 Choose the Riemannian metric provided by \cite[Lemma 10.5]{kirwan1985cohomology}, so that near each critical subset $C_\alpha$ of $f$, the gradient flow is tangent to its minimizing manifold $\Sigma_\alpha$. Although  \cite[Theorem 10.4]{kirwan1985cohomology} is stated for compact manifolds, Kirwan explains in  \cite[Section 9]{kirwan1985cohomology} that it remains valid in the noncompact case provided that every trajectory of the negative  gradient is contained in a compact subset. This holds here: if $\varphi_t(x)$ is the negative gradient trajectory through $x$ and $m$ is a lower bound for $f$, then $m\leq f\big(\varphi_t(x) \big)\leq f(x)$, so $\varphi_t(x) \subset f^{-1}\big([m,f(x)]  \big)$, which is compact since $f$ is proper. 

We apply Kirwan's Morse-stratification. Let $S_\alpha$ be the Morse stratum associated to $C_\alpha$. By  \cite[Theorem 10.4]{kirwan1985cohomology}, we have a smooth stratification $M=\sqcup_{\alpha}S_\alpha$, and $S_\alpha$ agrees with the minimizing manifold $\Sigma_\alpha$ near  $C_\alpha$.  Hence $\text{codim}(S_\alpha)=\text{codim}(\Sigma_\alpha)=\text{index}(C_\alpha)$.

Let $S_0 := \bigsqcup_{\operatorname{ind}(C_\alpha)=0} S_\alpha.$ Since $f$ has no points of index 1, we have $S_0=M\setminus \bigsqcup_{\text{ind}(C_\alpha)\geq 2}S_\alpha$. Since the complement of a locally finite union of smooth submanifolds of codimension at least  2 in a connected manifold is connected, $S_0$ is connected. Since  the Morse strata are pairwise disjoint and $S_0$ is connected, there is exactly one index 0 stratum, say $S$, with corresponding critical subset $C$. Thus $C_0=C$. By \cite[Lemma 10.17]{kirwan1985cohomology}, the inclusion $C \hookrightarrow S$ is an equivalence in \v{C}ech cohomology. Since $S$ is connected, $\check{H}^0(S;\mathbb{Z}/2)\cong\mathbb{Z}/2$. Therefore, $\check{H}^0(C;\mathbb{Z}/2) \cong
\check{H}^0(S;\mathbb{Z}/2) \cong \mathbb{Z}/2$, so $C=C_0$ is connected. \end{proof}

\section{Local Analysis near Elliptic-Regular Singularities}\label{sec:ER}
Let $\Phi:M\rightarrow \mathbb{R}^2$ be an almost-toric system on $(M,\omega)$, and let $f:\mathbb{R}^2 \rightarrow \mathbb{R}$ be a smooth function. In Section 3--5, we study the composition $H:=f\circ \Phi$ near each type of singularity of $\Phi$. Our goal is to establish conditions under which $H$ is Morse--Kirwan with no critical points of index 1, so that we can apply Lemma \ref{main lemma}. 

We begin with elliptic-regular singularities.  We first describe the local geometry of the rank 1 critical set and its image, which will be used in our analysis of $H$.

\begin{proposition}

Let $\Phi=(\Phi_1,\Phi_2) : M \to \mathbb{R}^2$ be a nondegenerate integrable system on $(M,\omega)$, and let  $C\subset \operatorname{Crit}_1(\Phi)$ be a connected component. Then $C$ is locally a 2-dimensional submanifold of $M$. Moreover, for every $p\in C$, there exists a neighborhood $U$ of $p$ such that $\Phi(C\cap U)$ is a smooth embedded curve in $\mathbb{R}^2$.  
\end{proposition}

\begin{proof}
 Let $p\in C$. Since  $\text{rank}(d\Phi)_p=1$, we may assume that $d\Phi_1(p)\neq 0$.
 By the Darboux-Carathéodory theorem (\cite[Theorem 5.2.19]{weinsteinr}, \cite[Theorem 4.1]{pelayo2011symplectic}), there exist local symplectic coordinates $x=(x_1,\xi_1,x_2,\xi_2)$ centered at $p$ on a sufficiently small neighborhood $U$ and $\xi_1=\Phi_1-\Phi_1(p)$. Since $\{\Phi_1,\Phi_2\}=0$, $\Phi_2$ is independent of  $x_1$. Moreover $\operatorname{rank}( d\Phi)_p=1$ implies 
 $\frac{\partial \Phi_2}{\partial x_2}(p)= \frac{\partial \Phi_2}{\partial \xi_2}(p)=0.$ Thus $\text{Crit}_1(\Phi)$ is locally given by $\text{Crit}_1(\Phi)=\{x \ | \ \frac{\partial \Phi_2}{\partial x_2}(x)= \frac{\partial \Phi_2}{\partial \xi_2}(x)=0 \}$.

 Since $p$ is nondegenerate, the matrix $\begin{pmatrix}
\frac{\partial^2 \Phi_2}{\partial x_2^2}(p) &
\frac{\partial^2 \Phi_2}{\partial x_2\partial \xi_2}(p)\\
\frac{\partial^2 \Phi_2}{\partial \xi_2\partial x_2}(p) &
\frac{\partial^2 \Phi_2}{\partial \xi_2^2}(p)
\end{pmatrix}$ is invertible. Hence, by the implicit function theorem, there exist unique smooth functions $x_2(\xi_1)$ and $\xi_2(\xi_1)$ defined near $\xi_1(p)=0$ such that $C\cap U=\{(x_1,\xi_1,x_2(\xi_1),\xi_2(\xi_1))\}$ for a sufficiently small neighborhood $U$ of $p$. Thus $C$ is locally a smooth two-dimensional submanifold of $M$.

Define $\tilde{\Phi}_2(\xi_1):=\Phi_2(\xi_1,x_2(\xi_1),\xi_2(\xi_1))$. Then $\Phi(C\cap U)$ is parametrized by $\xi_1 \mapsto (\xi_1+\Phi_1(p),\tilde{\Phi}_2(\xi_1)).$ Since the first component is $\xi_1+\Phi_1(p)$, this parametrization is an injective immersion, and its inverse onto its image is given by $(x,y)\mapsto x-\Phi_1(p)$. Therefore, $\Phi(C\cap U)$ is a smooth embedded curve in $\mathbb{R}^2$. \end{proof}

We next introduce some terminology needed to study the behavior of \(f\) near the elliptic-regular critical value curve of $\Phi$.

\begin{definition}\label{definition of restricted critical point}
    Let $\Phi:M\to\mathbb{R}^2$ be a nondegenerate integrable system on $(M, \omega)$, let $f:\mathbb{R}^2\rightarrow \mathbb{R}$ be a smooth function, and let $C\subset \operatorname{Crit}_1(\Phi)$ be a connected component. Fix $p\in C$, and choose a neighborhood $U$ of $p$ such that  $\Phi(C \cap  U)$ is a smooth embedded curve.   We say that $\Phi(p)$ is a \emph{critical point} of $f|_{\Phi(C\cap U)}$ if, for any smooth regular local parametrization $\delta:I\rightarrow \Phi(C\cap U)$ with $\delta(t_0)=\Phi(p)$, one has  $(f\circ \delta)'(t_0)=0$.
\end{definition}

\begin{definition}\label{definition of restricted MK}
   With the notation above, we say that the restriction $f|_{\Phi(C\cap U)}$ is a \emph{Morse-Kirwan function} near $\Phi(p)$ if, for any smooth regular local parametrization $\eta:I \rightarrow\Phi(C\cap U)$ with $\eta(t_0)=\Phi(p)$, the pullback $\eta^*(f|_{\Phi(C\cap U)})$  is a Morse-Kirwan function in a neighborhood of $t_0$.
\end{definition}

Since any two smooth regular local parametrizations of the same embedded curve differ by a local diffeomorphism, both definitions are independent of the choice of parametrization.

\vspace{1mm}

\begin{definition}
Let $p\in M$ be an elliptic-regular critical point of a nondegenerate integrable system $\Phi:M\to\mathbb R^2$. We say that a nonzero vector $v\in T_{\Phi(p)}\mathbb R^2$ \emph{points locally into} the interior of $\Phi(M)$ at $\Phi(p)$  if, for every sufficiently small neighborhood $U$ of $p$, there
exists  $\varepsilon>0$ such that
\[
\Phi(p)+tv \in \text{Int}\big(\Phi(U)\big) \qquad \text{for all } t\in(0,\varepsilon).
\] 
We say that \(v\) points \emph{locally outward} from \(\Phi(M)\) at \(\Phi(p)\) if it does not point locally into the interior of \(\Phi(M)\) at \(\Phi(p)\).

\end{definition}

\medskip

We now establish the local setup used throughout the remainder of this section. Let $p$ be an elliptic-regular critical point of a nondegenerate integrable system $\Phi:M\rightarrow \mathbb{R}^2$ on $(M,\omega)$. By Theorem \ref{local normal form theorem}, there exist local symplectic coordinates $x=(x_1,\xi_1,x_2,\xi_2)$ centered at $p$ on a sufficiently small neighborhood $U$, and a local diffeomorphism $g:\mathbb{R}^2\rightarrow \mathbb{R}^2$ with $g(0)=\Phi(p)$ such that  
\[\Phi=g\circ q \quad \text{where }\quad q=(q_1,q_2)=(\frac{x_1^2+\xi_1^2}{2}, \xi_2). \]
Let $f:\mathbb{R}^2\rightarrow \mathbb{R}$ be a smooth function, and consider a composite function $H:=f\circ\Phi=(f\circ g)(q_1,q_2)$ on $U$.
Then \[dH=\frac{\partial(f\circ g)}{\partial q_1}(x_1dx_1+\xi_1d\xi_1)+\frac{\partial(f\circ g)}{\partial q_2}d\xi_2.\]
Hence, if  $p\in \text{Crit}(H)$ and $\Phi(p)\notin \text{Crit}(f)$, then 
$\frac{\partial(f\circ g)}{\partial q_2}(0)=  0$ and $\lambda_1:=\frac{\partial(f\circ g)}{\partial q_1}(0)\neq0$. The following lemma relates the sign of $\lambda_1$ to the direction of the gradient vector $\nabla f_{\Phi(p)}$ of $f$ at $\Phi(p)$.

\medskip

\begin{lemma}\label{sign of lambda_1 and inward pointing}
With the notation above, suppose that $p\in \operatorname{Crit}(H)$ and $\Phi(p)\notin \operatorname{Crit}(f)$. Then $\lambda_1>0$ if and only if $\nabla f(\Phi(p))$ points locally inward into the interior of $\Phi(M)$, and $\lambda_1<0$ if and only if $\nabla f(\Phi(p))$ points locally outward from  $\Phi(M)$ .
\end{lemma}

\begin{proof}

Let $v_1:=dg_0(\partial_{q_1}),\  v_2:=dg_0(\partial_{q_2})$. Since $g$ is a local diffeomorphism, $v_1, v_2$  form a basis of $T_{\Phi(p)}\mathbb R^2$. 
Hence there exist unique $a,b\in\mathbb R$
such that $\nabla f_{\Phi(p)}=a v_1+b v_2.$
Since $p\in \text{Crit}(H)$, we have $0=\frac{\partial(f\circ g)}{\partial q_2}(0)=\langle \nabla f_{\Phi(p)},v_2\rangle$.
 This gives $0=\langle a v_1+b v_2,\,v_2\rangle
=a\langle v_1,v_2\rangle+b\langle v_2,v_2\rangle$, thus $b=-a\,\frac{\langle v_1,v_2\rangle}{\langle v_2,v_2\rangle}.$ Therefore  $\nabla f_{\Phi(p)}=a\Big(v_1-\frac{\langle v_1,v_2\rangle}{\langle v_2,v_2\rangle}v_2\Big)
=:a\,v_2^\perp$
where $v_2^\perp\perp v_2$. Consequently,
\[
\lambda_1=\frac{\partial(f\circ g)}{\partial q_1}(0)
=\langle \nabla f_{\Phi(p)},v_1\rangle
=a\Big(\langle v_1,v_1\rangle
-\frac{\langle v_1,v_2\rangle^2}{\langle v_2,v_2\rangle}\Big)=a||v_2^\perp||^2.
\]
Thus $\operatorname{sign}(\lambda_1)=\operatorname{sign}(a)$. Moreover, $a\neq 0$ since $\Phi(p)\notin \text{Crit}(f)$.

Now define $\gamma(t):=g^{-1}\big(\Phi(p)+t\nabla f_{\Phi(p)}\big)$ for sufficiently small $t\in \mathbb{R}$. Since $g(0)=\Phi(p)$, we have $\gamma(0)=0$. Since $d(g^{-1})_{\Phi(p)}(v_i)=\partial_{q_i}$ with $i=1,2$, we have   $\gamma'(0)=d(g^{-1})_{\Phi(p)}(\nabla f_{\Phi(p)})=a\partial_{q_1}+b\partial_{q_2}$. Taylor expansion at $t=0$ gives $\gamma(t)=(at, bt)+O(t^2)$. In particular, the first component of $\gamma(t)$ is $\gamma_1(t)= at+O(t^2)$.  

For every sufficiently small neighborhood $U$ of $p$, the elliptic-regular local normal form implies that, near $\Phi(p)$, the interior of $\Phi(U)$ corresponds under $g^{-1}$ to $\{(q_1,q_2) \; | \; q_1>0\}$ near the origin. Hence, for sufficiently small $t>0$, $\nabla f _{\Phi(p)}$ points locally inward into the interior of $\Phi(M)$ if and only if $a>0$, and locally outward if and only if $a<0$. Since $\operatorname{sign}(\lambda_1)=\operatorname{sign}(a)$, the result follows. 
\end{proof}

Now we give sufficient conditions for $H=f\circ \Phi$ to be Morse--Kirwan near an elliptic-regular critical point and determine its possible indices.

\begin{theorem}\label{ER MK}
Let $\Phi:M\to\mathbb{R}^2$ be a nondegenerate integrable system on $(M,\omega)$, and let $f:\mathbb{R}^2\to\mathbb{R}$ be a smooth function.  Let $p$ be an elliptic-regular critical point of $\Phi$, and let $C\subset \operatorname{Crit}_1(\Phi)$ be a connected component containing $p$. Suppose: 

\begin{itemize}
\item $\Phi(p)\notin \operatorname{Crit}(f)$.
\item The restriction $f|_{\Phi(C)}$ is a Morse--Kirwan function near $\Phi(p)$.

\end{itemize}

Then $H:=f\circ \Phi$ is a Morse--Kirwan function  near $p$. Moreover, if $p \in \operatorname{Crit}(H)$, 
the Morse index of $H$ at $p$ is determined as follows: 
\[
\operatorname{index}_p(H) =
\begin{cases}
0 & \text{if } \Phi(p) \text{ is a local minimum of $f|_{\Phi(C)}$ and } \nabla f_{\Phi(p)} \text{ points locally inward into} \ \Phi(M), \\
1 & \text{if } \Phi(p) \text{ is a local maximum of $f|_{\Phi(C)}$ and } \nabla f_{\Phi(p)} \text{ points locally inward into} \ \Phi(M),\\
2 & \text{if } \Phi(p) \text{ is a local minimum of $f|_{\Phi(C)}$ and } \nabla f_{\Phi(p)} \text{ points locally outward from} \ \Phi(M), \\
3 & \text{if } \Phi(p) \text{ is a local maximum of $f|_{\Phi(C)}$ and } \nabla f_{\Phi(p)} \text{ points locally outward from} \ \Phi(M).
\end{cases}
\]
\end{theorem}

\begin{proof}
We work in the local coordinates and notation established above. Recall that $H=f\circ\Phi=(f\circ g)(q_1,q_2)$ where $q_1=\frac{x_1^2+\xi_1^2}{2}, \ q_2=\xi_2$. Since $dH=\frac{\partial(f\circ g)}{\partial q_1}(x_1dx_1+\xi_1d\xi_1)+\frac{\partial(f\circ g)}{\partial q_2}d\xi_2$, we have $\text{Crit}(H)\cap U=\Big\{(x_1,\xi_1,x_2,\xi_2)\in U: \  x_1\frac{\partial (f\circ g)}{\partial q_1}=\xi_1\frac{\partial (f\circ g)}{\partial q_1}=\frac{\partial (f\circ g)}{\partial q_2}=0\Big\}.$ Since $\Phi(p)\notin \operatorname{Crit}(f)$, after shrinking $U$ if necessary, we may assume  $d(f \circ g)=\bigg[\frac{\partial (f\circ g)}{\partial q_1}\ \ \frac{\partial (f\circ g)}{\partial q_2}\bigg] \neq 0$ on $q(U)$. Hence
\[ 
\operatorname{Crit}(H)\cap U = \{(0,0,x_2,\xi_2)\in U:\frac{\partial (f\circ g)}{\partial q_2}(0,q_2)=0 \}.\] 


Define a function $\varphi(q_2):=(f\circ g)(0,q_2)$. After shrinking $U$ if necessary, the map $q_2\mapsto g(0,q_2)$ is a smooth regular local parametrization of $\Phi(C\cap U)$ near $\Phi(p)$. Hence, by the hypothesis that $f|_{\Phi(C)}$ is Morse-Kirwan near $\Phi(p)$,  $\varphi$  is Morse-Kirwan near 0. Thus, after restricting to a neighborhood of $0$, its critical set is a locally finite disjoint union  $\text{Crit}(\varphi)=\bigsqcup_{\alpha}C_{\alpha}$ where each $C_\alpha$ is a closed subset on which $\varphi$ is constant and admits a minimizing manifold  $N_\alpha$ for $\varphi$ along $C_\alpha$. Then  the corresponding critical subset of $H$ is
\[\operatorname{Crit}(H)\cap U=\bigsqcup_{\alpha}\widetilde{C}_\alpha \quad \text{where} \quad \widetilde{C}_\alpha:=\left\{(0,0,x_2,\xi_2)\in U:\;\xi_2\in C_\alpha\right\}.\]
This decomposition is locally finite, and $H$ is constant on each $\widetilde{C}_\alpha$.

Fix $\alpha$. Since $d(f\circ g)\neq 0$ on $q(U)$, after shrinking $U$ we may assume that $\lambda_1(q_2):= \frac{\partial (f \circ g)}{\partial q_1}(0,q_2)$ has constant nonzero sign. At each  $x=(0,0,x_2,\xi_2)\in \widetilde{C}_\alpha$, we have 
\[\text{Hess}_x(H)=\lambda_1(\xi_2)\big(dx_1\otimes dx_1+d\xi_1\otimes d\xi_1\big)+\varphi''(\xi_2)\;d\xi_2\otimes d\xi_2\] and $\partial_{x_2}$ lies in its kernel. The eigenvalues of $\text{Hess}_x(H)$ are given by $\{\lambda_1(\xi_2),\lambda_1(\xi_2),\varphi''(\xi_2),0\}$.

\vspace{2mm}

Now we construct a minimizing manifold $\widetilde{N}_\alpha$ for $H$ along $\widetilde{C}_\alpha$.
\vspace{2mm}






\noindent\emph{Case 1: $\lambda_1>0$.}

\noindent Define $\widetilde{N}_\alpha:=\{(x_1,\xi_1, x_2, \xi_2)\in U \; | \;  \xi_2\in N_\alpha \; \}.$ Since $\lambda_1>0$, $\text{Hess}_x(H)$ is positive definite on $\operatorname{span}\{\partial_{x_1},\partial_{\xi_1}\}$. Moreover since $T_{\xi_2}N_\alpha$ is maximal among subspaces on which $\text{Hess}_{\xi_2}(\varphi)$ is positive semidefinite, it follows that $T_x\widetilde{N}_\alpha$ is maximal  subspaces on which $\text{Hess}_x$(H) is positive semidefinite.  Moreover, since $q_1\geq 0$ and $\lambda_1>0$ near the origin,  we have
\[H=(f\circ g)(q_1,\xi_2)\ge (f\circ g)(0,\xi_2)=\varphi(\xi_2). \]
Since $\varphi|_{N_\alpha}$ attains its minimum exactly on $C_\alpha$, it follows that $H|_{\widetilde{N}_\alpha}$ attains its minimum exactly on $\widetilde{C}_\alpha$. Thus $\widetilde{N}_\alpha$ is a minimizing manifold for $H$ along $\widetilde{C}_\alpha$.
\vspace{2mm}

\noindent\emph{Case 2: $\lambda_1<0$.}

\noindent Define $\widetilde{N}_\alpha:=\{(0,0, x_2, \xi_2)\in U \; | \;  \xi_2\in N_\alpha \; \}.$ Since $\lambda_1<0$, 
$\text{Hess}_x(H)$ is negative definite on $\operatorname{span}\{\partial_{x_1},\partial_{\xi_1}\}$. Moreover since  $T_{\xi_2}N_\alpha$ is maximal among subspaces on which $\text{Hess}_{\xi_2}(\varphi)$ is positive semidefinite, it follows that $T_x\widetilde{N}_\alpha$ is maximal subspaces on which $\text{Hess}_x$(H) is positive semidefinite.  Moreover we have 
\[H|_{\widetilde{N}_\alpha}=\varphi(\xi_2)|_{N_\alpha}, \]
which attains its minimum exactly on $\widetilde{C}_\alpha$. Thus $\widetilde{N}_\alpha$ is a minimizing manifold for $H$ along $\widetilde{C}_\alpha$. Therefore $H$ is Morse--Kirwan near $p$.

Lastly, suppose $p\in \text{Crit}(H)$. The above constructions implies 
\[
\operatorname{index}_p(H)
=
\begin{cases}
\operatorname{index}_0(\varphi),&\lambda_1>0,\\[2mm]
\operatorname{index}_0(\varphi)+2,&\lambda_1<0.
\end{cases}
\]
Since $\varphi$ is a one-dimensional Morse--Kirwan function, its index at 0 is either 1 or 1. Index 0 corresponds to a local minimum (possibly degenerate), while index 1 corresponds to a nondegenerate local maximum.  By Lemma \ref{sign of lambda_1 and inward pointing},  $\lambda_1>0$ if and only if $\nabla f_{\Phi(p)}$ point locally inward into $\Phi(M)$, whereas $\lambda_1<0$ if and only if $\nabla f_{\Phi(p)}$ points locally outward from $\Phi(M)$. This gives the stated index. 

\end{proof}



\begin{example}
Consider Example 4.10 of \cite{pelayo2015fiber}, given by
\[
M=\mathbb S^2\times \mathbb S^1\times \mathbb S^1,
\qquad
\Phi(h,a,b,c)=\bigl(h\cos(nb),\,h\sin(nb)\bigr),
\]
where \(h\in[1,2]\), \(a\in \mathbb R/2\pi\mathbb Z\) is the angular coordinate on
\(\mathbb S^2\),  \(b,c\in \mathbb R/2\pi\mathbb Z\), and $n$ is a positive integer.
Then
\[
\Phi(M)=\{(x,y)\in \mathbb R^2 \mid 1\le x^2+y^2\le 4\}.
\]
Define $f(x,y):=\left\|(x,y)-\left(\frac{3}{2},0\right)\right\|^2.$
Then the level set \(f^{-1}((\frac{5}{2})^2)\) is tangent to the inner boundary of \(\Phi(M)\) at $q=(-1,0).$  Along the inner boundary circle, \(q\) is the farthest point from the center \(\left(\frac{3}{2},0\right)\). Hence the restriction of \(f\) to the inner boundary has a \emph{local maximum} at \(q\). Moreover, $\nabla f(q)$ \emph{points locally inward} into \(\Phi(M)\). Thus Theorem \ref{ER MK} implies that $\text{index}_p(f\circ \Phi)=1 $ where $p$ is an elliptic-regular point such that $\Phi(p)=q$.

We can also verify this directly. 
Let $p=(\text{south pole},\,b_0,\,c_0)\in M$,  $nb_0\equiv \pi \pmod{2\pi}$ with arbitrary \(c_0\in \mathbb S^1\). Then \(\Phi(p)=q\). Writing $H:=f\circ \Phi$, we have 
\[H(h,a,b,c)=\left(h\cos(nb)-\frac32\right)^2+h^2\sin^2(nb).\] A direct calculation shows that  $\text{index}_p(f\circ \Phi)=1 $.
\end{example}

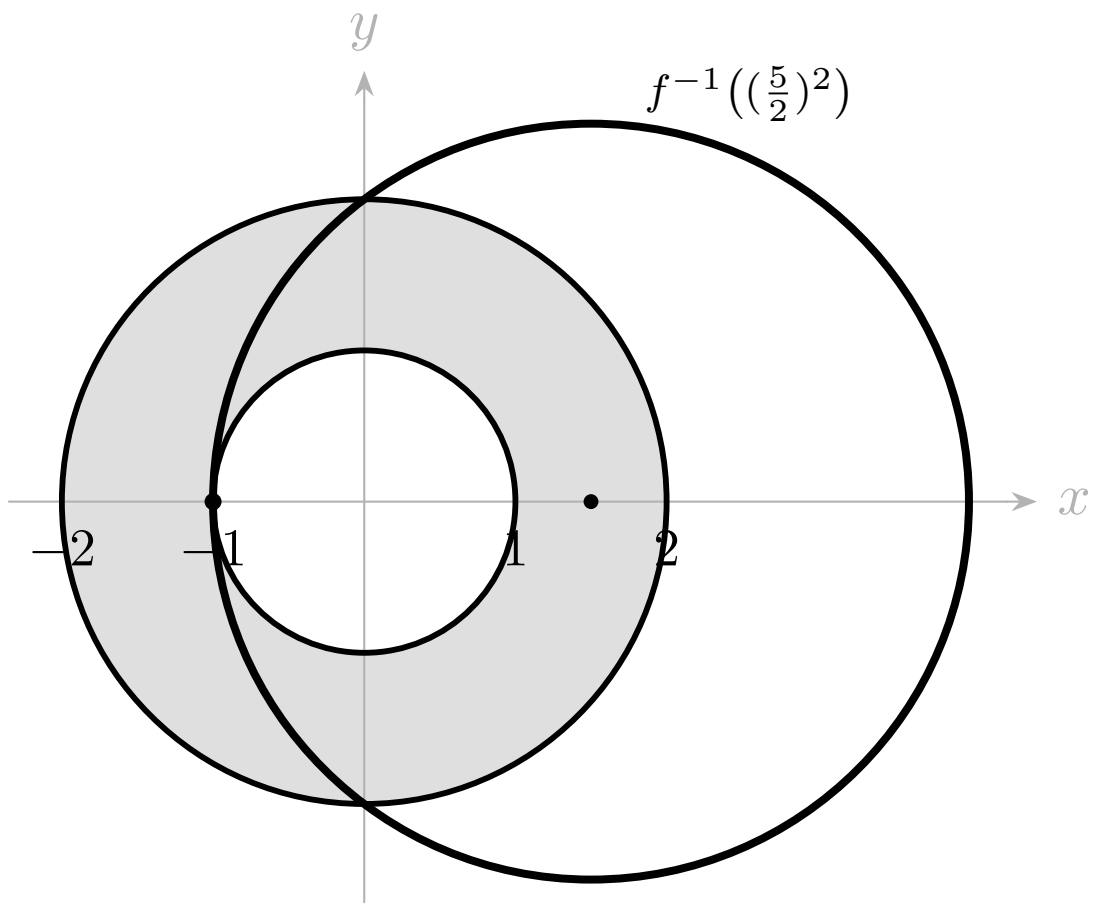
\begin{figure}[ht]
\centering
\begin{tikzpicture}[scale=1.00, line cap=round, line join=round, >=Stealth]

  \fill[gray!25, even odd rule] (0,0) circle (2) (0,0) circle (1);

  \draw[gray!60] (-2.35,0) -- (4.25,0);
  \draw[gray!60] (0,-2.65) -- (0,2.65);
  \draw[->,gray!60] (4.25,0) -- (4.45,0) node[right] {$x$};
  \draw[->,gray!60] (0,2.65) -- (0,2.85) node[above] {$y$};

  \foreach \x in {-2,-1,1,2}{
    \draw[gray!60] (\x,0.06) -- (\x,-0.06);
    \node[below] at (\x,-0.06) {\small $\x$};
  }

  \draw[line width=1.1pt] (0,0) circle (2);
  \draw[line width=1.1pt] (0,0) circle (1);

  \draw[line width=1.5pt] (1.5,0) circle (2.5);

  \fill (1.5,0) circle (1.4pt);

  \fill (-1,0) circle (1.6pt);

  \node at (2.55,2.70) {\scriptsize $f^{-1}\big((\frac{5}{2})^2\big)$};

\end{tikzpicture}
\caption{The image \(\Phi(M)\) in Example 4.10 of \cite{pelayo2015fiber} and the level set \(f^{-1}\big((\frac{5}{2})^2\big)\), which is tangent to the inner boundary at \((-1,0)\).}
\label{fig:index1-annulus}
\end{figure}

Next, we give a sufficient condition for $H=f\circ\Phi$ to be Morse–Kirwan near an elliptic-regular critical point $p$, allowing $\Phi(p)$ to be a critical point of $f$.

\begin{theorem}\label{ER MK local min}
Let $\Phi:M\to\mathbb{R}^2$ be a nondegenerate integrable system on $(M,\omega)$, and let
$f:\mathbb{R}^2\to\mathbb{R}$ be a Morse--Kirwan function.  Let  $p$ be an elliptic-regular critical point of $\Phi$.
Suppose:

\begin{itemize}
\item  $\Phi(p)$ is a local minimum of $f$.
\item The restriction $f|_{\Phi(C)}$ is a Morse-Kirwan function near $\Phi(p)$ where $C \subset \operatorname{Crit}_1(\Phi)$ is a connected component containing $p$.
\end{itemize}

Then $f\circ\Phi$ is a Morse--Kirwan function near $p$, with index 0 at every critical point sufficiently near $p$.
\end{theorem}

\begin{proof} Let $H := f \circ \Phi $ and we may assume $H(p) = f(\Phi(p)) = 0$. Since $\Phi(p)$ is a local minimum of $f$, continuity of $\Phi$ implies that  $p$ is a local minimum of $H$. Working in the local coordinates introduced above on a neighborhood $U$ of $p$, we may assume, after shrinking $U$ if necessary, that $H\geq 0$ on $U$. Moreover we have $H=f\circ\Phi=(f\circ g)(q_1,q_2)$, and 
\[ dH = \frac{\partial (f\circ g)}{\partial q_1}(x_1\,dx_1 + \xi_1\,d\xi_1) + \frac{\partial (f\circ g)}{\partial q_2}d\xi_2.
\]
Define $\varphi(q_2):=(f\circ g)(0,q_2)$. Then, after shrinking $U$ if necessary, we have 
\[
\mathrm{Crit}(H)\cap U=A\cup B
\]
where $A:=\Phi^{-1}(\mathrm{Crit}(f))\cap U$ and $B:=\{(0,0,x_2,\xi_2)\in U\; | \;  \varphi'(\xi_2)=0\}$.

\vspace{2mm}

We claim that $A\cup B \subset H^{-1}(0)$. 
We first show  $A\subset H^{-1}(0)$. Since $f$ is Morse–Kirwan and $\Phi(p)$ is a critical point of $f$ with $f(\Phi(p))=0$, Proposition \ref{equivalence of locally finitness} gives a neighborhood $V$ of $\Phi(p)$ such that $\mathrm{Crit}(f)\cap V \subset f^{-1}(0).$ After shrinking $U$ so that $\Phi(U)\subset V$, we obtain
\[
A=\Phi^{-1}(\mathrm{Crit}(f))\cap U \subset
\Phi^{-1}(f^{-1}(0))
=
H^{-1}(0).
\] 

Now we show $B\subset H^{-1}(0)$. Since the map $q_2\mapsto g(0,q_2)$ parametrizes $\Phi(C)$ near $\Phi(p)$, the hypothesis that $f|_{\Phi(C)}$ is Morse-Kirwan  near $\Phi(p)$ implies that $\varphi$ is Morse-Kirwan near 0. Moreover, 
since $\Phi(p)$ is a local minimum of $f$ and
$\varphi(0)=(f\circ g)(0,0)=f(\Phi(p))=0$, 0 is a local minimum of $\varphi$. Applying Proposition \ref{equivalence of locally finitness} to $\varphi$, there exists a neighborhood $W$ of 0 in $q_2(U)$ such that  $\text{Crit}(\varphi)\cap W \subset \varphi^{-1}(0)$. After shrinking $U$ so that $q_2(U)\subset W$, every $x\in B$ satisfies
\[H(x)=f(g(0,\xi_2))=\varphi(\xi_2)=0.\]
Hence $B\subset H^{-1}(0).$ It follows that $\text{Crit}(H)\cap U=A\cup B\subset H^{-1}(0)$.

Since  $p$ is a local minimum of $H$ and $\text{Crit(H)}\cap U \subset H^{-1}(0)$, Lemma \ref{local min and preimage} implies that $H$ is Morse-Kirwan near $p$. Moreover, every $x\in \text{Crit}(H)\cap U$ is a local minimum of $H$ since  $\text{Crit}(H)\cap U\subset H^{-1}(0)$ and $H\geq0$ on $U$.  Thus $\text{Hess}_x(H)$ is positive semidefinite at every $x\in \text{Crit}(H)\cap U$, and hence $H$ has index $0$ on $\text{Crit}(H)\cap U$. \end{proof}

\section{Local Analysis near Elliptic-Elliptic Singularities}\label{sec:EE}

In this section, we study the composition $H=f\circ\Phi$ near elliptic-elliptic singularities. We first describe the local geometry of the rank 1 critical set near an elliptic-elliptic critical point.

\begin{proposition}
Let $\Phi:M\to\mathbb{R}^2$ be an almost-toric system on  $(M,\omega)$. Let $C$ be a connected component of $\operatorname{Crit}_1(\Phi)$. If $p\in \overline{C}\setminus C$, then $p$ is an elliptic-elliptic critical point of $\Phi$.
\end{proposition}

\begin{proof}
Since $p\in \overline{C}$, there exists a sequence $p_n\in C$ such that $p_n\rightarrow p$. Since the set of regular points of $\Phi$ is open, $p$ cannot be regular. Suppose $\text{rank}(d\Phi_p)=1$.  Then $C\cup \{p\}$ is a connected subset of $\text{Crit}_1(\Phi)$. Since $C$ is a connected component of $\text{Crit}_1(\Phi)$, it follows that $p\in C$, contradicting $p\in \overline{C}\setminus C$.  Hence  $\text{rank}(d\Phi_p)=0$. Since focus-focus critical points are  isolated in $\text{Crit}(\Phi)$, $p$ cannot be focus-focus. Therefore $p$ is elliptic-elliptic.

\end{proof}

\begin{proposition}\label{smoothness at endpoint EE}
Let $\Phi : M \to \mathbb{R}^2$ be an almost-toric system on $(M,\omega)$, and let $p$ be an elliptic--elliptic critical point of $\Phi$. Then there exists a neighborhood $U$ of $p$ such that $\operatorname{Crit}_1(\Phi) \cap U$ has exactly two connected components $C_1, C_2$, whose closures in $U$ contain $p$. Moreover, for $i=1,2$,  $\Phi(\overline{C_i})$  is locally a smooth embedded curve  near $\Phi(p)$, with boundary point $\Phi(p)$.
\end{proposition}

\begin{proof}
   By Theorem \ref{local normal form theorem}, there exist local symplectic coordinates $(x_1,\xi_1,x_2,\xi_2)$ centered at $p$ on a sufficiently small neighborhood $U$, and a local diffeomorphism $g : \mathbb{R}^2 \to \mathbb{R}^2$ with $g(0) = \Phi(p)$, such that $\Phi = g \circ q$ where $q = (q_1, q_2) =(\frac{x_1^2 + \xi_1^2}{2},  \frac{x_2^2 + \xi_2^2}{2}).$ Since $g$ is a local diffeomorphism, $\text{rank}(d\Phi)=\text{rank}(dq)$. Thus 
$\text{Crit}_1(\Phi)\cap U=C_1\sqcup C_2$ where 
\[C_1:=\{(0,0,x_2,\xi_2)\in U \; | \; (x_2,\xi_2)\neq(0,0)\}, \quad C_2:=\{(x_1,\xi_1,0,0)\in U \; | \; (x_1,\xi_1)\neq(0,0)\}.\]
After choosing $U$ sufficiently small, both $C_1$ and $C_2$ are connected, and their closures in $U$ are 
\[
 \overline{C_1}=\{(0,0,x_2,\xi_2)\in U\}
\quad \text{and} \quad
\overline{C_2}=\{(x_1,\xi_1,0,0)\in U\},
\]
both of which contain $p$.  Moreover, near the origin, we have $q(\overline{C_1})=\{(0,s)\;|\; s\geq0 \}$ and  $q(\overline{C_2})=\{(s,0)\;|\; s\geq0 \}$. Since $g$ is a local diffeomorphism, $\Phi(\overline{C_i})=g\big(q(\overline{C_i})\big)$  is locally a smooth embedded curve with boundary  $g(0)=\Phi(p)$ for $i=1,2$.\end{proof}

By Proposition \ref{smoothness at endpoint EE}, Definition \ref{definition of restricted critical point} and Definition \ref{definition of restricted MK}  extend naturally to $f|_{\Phi(\overline{C_i})}$ near the elliptic-elliptic critical value $\Phi(p)$, using smooth regular local parametrizations $\eta:[0,\epsilon)\rightarrow \Phi(\overline{C_i})$ with $\eta(0)=\Phi(p)$.

We now give sufficient conditions for $H=f\circ\Phi$ to be Morse–Kirwan near an elliptic-elliptic critical point $p$ when $\Phi(p)$ is not a critical point of $f$.

\begin{theorem}\label{EE MK}
Let $\Phi:M\to\mathbb{R}^2$ be an almost-toric system on  $(M,\omega)$, and let $f:\mathbb{R}^2\to\mathbb{R}$ be a smooth function. Let $p$ be an elliptic-elliptic critical point of $\Phi$, and suppose that $\Phi(p)\notin \operatorname{Crit}(f)$. For each connected component $C\subset \operatorname{Crit}_1(\Phi)$ such that $p\in \overline{C}$, suppose:

\begin{itemize}
\item The restriction $f|_{\Phi(\overline{C})}$ is a Morse-Kirwan function near $\Phi(p)$.
\item  $\Phi(p)$ is either a regular point or  a local minimum of  $f|_{\Phi(\overline{C})}$.

\end{itemize}

Then $H:=f\circ\Phi$ is a Morse-Kirwan function  near $p$ with index  0, 2, or 4 at $p$.
\end{theorem}

\begin{proof}
 By Theorem \ref{local normal form theorem}, there exist local symplectic coordinates $(x_1,\xi_1,x_2,\xi_2)$ centered at $p$ on a sufficiently small neighborhood $U$, and a local diffeomorphism $g : \mathbb{R}^2 \to \mathbb{R}^2$ with $g(0) = \Phi(p)$, such that $\Phi = g \circ q$ where $q = (q_1, q_2) =(\frac{x_1^2 + \xi_1^2}{2},  \frac{x_2^2 + \xi_2^2}{2}).$ Let $H:=f\circ\Phi=(f\circ g)(q_1,q_2)$. Then  \[dH=\frac{\partial(f\circ g)}{\partial q_1}(x_1dx_1+\xi_1d\xi_1)+\frac{\partial(f\circ g)}{\partial q_2}(x_2dx_2+\xi_2d\xi_2).\] Since 
$\Phi(p)\notin \text{Crit}(f)$ and $g$ is a local diffeomorphism $d(f\circ g)_0=\bigg[\frac{\partial (f\circ g)}{\partial q_1}(0)\ \ \frac{\partial (f\circ g)}{\partial q_2}(0)\bigg] \neq0$. Hence at least one of $\frac{\partial (f\circ g)}{\partial q_1}(0), \frac{\partial (f\circ g)}{\partial q_2}(0)$
is nonzero.

\vspace{3mm}

\noindent\emph{Case 1:} Suppose  that both $\frac{\partial (f\circ g)}{\partial q_1}(0)$ and $\frac{\partial (f\circ g)}{\partial q_2}(0)$ are nonzero. After shrinking $U$ if necessary, both remain nonzero on $q(U)$. Thus $\text{Crit}(H)\cap U=\{p\}$. Moreover, 
\[\text{Hess}_p(H)=\lambda_1(dx_1\otimes dx_1+d\xi_1\otimes d\xi_1)+\lambda_2(dx_2\otimes dx_2+d\xi_2\otimes d\xi_2),\]
where $\lambda_i=\frac{\partial (f\circ g)}{\partial q_i}(0)$, $i=1,2$. Thus its eigenvalues are $\{\lambda_1, \lambda_1,\lambda_2,\lambda_2\}$, hence  $p$ is a nondegenerate critical point of index 0, 2, or 4, depending on the sign of $\lambda_1,\lambda_2$. Hence $H$ is Morse, and therefore Morse-Kirwan, near $p$.

\vspace{3mm}
\noindent\emph{Case 2:}
Suppose that one of $\frac{\partial (f\circ g)}{\partial q_1}(0)$, $\frac{\partial (f\circ g)}{\partial q_2}(0)$  is nonzero and the other is zero. Without loss of generality, we may assume  $\frac{\partial (f\circ g)}{\partial q_1}(0)\neq 0$ and $\frac{\partial (f\circ g)}{\partial q_2}(0)=0$. After shrinking $U$ if necessary, we may assume that $\frac{\partial (f\circ g)}{\partial q_1}$ has constant nonzero sign on $q(U)$.

Define $\varphi(q_2):=(f\circ g)(0,q_2)$, $q_2\geq 0$. The map $q_2\mapsto g(0,q_2)$  parametrizes $\Phi(\overline{C_1})$ near $\Phi(p)$ where $C_1:=\{(0,0,x_2,\xi_2)\in U \; | \;(x_2,\xi_2)\neq(0,0)\} \subset \operatorname{Crit}_1(\Phi)$. Therefore the hypothesis that $f|_{\Phi(\overline{C_1})}$ is Morse-Kirwan  near $\Phi(p)$ implies that  $\varphi$ is Morse-Kirwan near $0$. Thus, after restricting to a neighborhood of $0$, its critical set is a locally finite disjoint union  $\text{Crit}(\varphi)=\bigsqcup_{\alpha}C_{\alpha}$ where each $C_\alpha$ is a closed subset on which $\varphi$ is constant.  Since  $\varphi'(0)=\frac{\partial (f\circ g)}{\partial q_2}(0)=0$, the point $\Phi(p)$ is not a regular point of $f|_{\Phi(\overline{C_1})}$. By hypothesis, $\Phi(p)$ is therefore a local minimum, and hence $0$ is a local minimum of $\varphi$. Thus we may assume that $\varphi\geq \varphi(0)$ near 0.

Since the family $\{C_\alpha\}$ is locally finite and consists of pairwise disjoint closed subsets, after shrinking $U$ if necessary, we may assume that $\text{Crit}(\varphi)\cap q_2(U)= C_\alpha\cap q_2(U)$, where $C_\alpha$ contains 0. Moreover, since $\varphi$ is Morse--Kirwan, Proposition \ref{equivalence of locally finitness} implies $\text{Crit}(\varphi)\cap q_2(U)\subset \varphi^{-1}(\varphi(0))$. Thus every critical point of $\varphi$ in $q_2(U)$ is a \emph{local minimum}, so $\varphi''(q_2)\geq 0$ on  $\text{Crit}(\varphi)\cap q_2(U)$.

Since $\frac{\partial (f\circ g)}{\partial q_1}\neq 0$, $dH=0$ implies  
\[\text{Crit}(H)\cap U= \widetilde{C}_\alpha \quad \text{where} \quad \widetilde{C}_\alpha:=\{(0,0,x_2,\xi_2)\in U \;|\; q_2\in C_\alpha\}.\] The set $\widetilde{C}_\alpha$ is closed in $U$, and $H$ is constant on it. 
At each $x\in \widetilde{C}_\alpha$, we have  \[\text{Hess}_x(H)=\lambda_1(q_2) \big(dx_1\otimes dx_1 +d\xi_1\otimes d\xi_1\big)+\varphi''(q_2)\;dq_2\otimes dq_2\] where $\lambda_1(q_2)=\frac{\partial (f\circ g)}{\partial q_1}(0,q_2)$,  and its eigenvalues are $\{\lambda_1(q_2), \; \lambda_1(q_2), \; 2q_2\varphi''(q_2), \; 0\}$. Since $q_2\geq 0$ and $\varphi''(q_2)\geq0$  on $C_\alpha\cap q_2(U)$, the eigenvalue $2q_2\varphi''(q_2)$ is nonnegative. In particular, at $p$, the eigenvalues are $\{\lambda_1(0), \lambda_1(0), 0,0\}$.

\vspace{2mm}

Now we construct a minimizing manifold $\widetilde{N}_\alpha$ for $H$ along $\widetilde{C}_\alpha$.

\vspace{1mm}

\noindent\emph{Case 2-(a): $\lambda_1>0$.}

\noindent Define $\widetilde{N}_\alpha:=U$. Since $\lambda_1>0$ and  $\varphi''(q_2)\geq 0 $ on $C_\alpha\cap q_2(U)$, $\text{Hess}_x(H)$ is positive semidefinite on $T_x\widetilde{N}_\alpha$ for every $x\in \widetilde{C}_\alpha$. Thus   $T_x\widetilde{N}_\alpha$ is a maximal among subspaces on which $\text{Hess}_x(H)$ is positive semidefinite. Moreover, since  $q_1\geq0$ and $\lambda_1>0$ on $q(U)$, 
\[H=(f\circ g)(q_1,q_2)\geq (f\circ g)(0,q_2)=\varphi(q_2)\geq \varphi|_{C_\alpha}=\varphi(0).\] Thus $H|_{\widetilde{N}_\alpha}=H$ attains its minimum exactly on $\widetilde{C}_\alpha$. Hence $\widetilde{N}_\alpha$ is a minimizing manifold for $H$ along $\widetilde{C}_\alpha$.

\vspace{1mm}

\noindent\emph{Case 2-(b): $\lambda_1<0$.}

\noindent Define $\widetilde{N}_\alpha:=\{(0,0,x_2,\xi_2)\in U\}$.  Since $\lambda_1<0$, $\text{Hess}_x(H)$ is negative definite on $\text{span}\{\partial_{x_1},\partial_{\xi_1}\}$ for every $x\in \widetilde{C}_\alpha$.
Morevover, since  $\varphi''(q_2)\geq 0$ on $C_\alpha\cap q_2(U)$, $\text{Hess}_x(H)$ is positive semidefinite on $T_x\widetilde{N}_\alpha=\text{span}\{\partial_{x_2},\partial_{\xi_2}\}$. Thus $T_x\widetilde{N}_\alpha$ is  maximal among subspaces in which $\text{Hess}_x(H)$ positive semidefinite. Furthermore, $H|_{\widetilde{N}_\alpha}=\varphi(q_2)\geq \varphi|_{C_\alpha}= \varphi(0)$, which attains its minimum exactly on $\widetilde{C}_\alpha$. Hence $\widetilde{N}_\alpha$ is a minimizing manifold for $H$ along $\widetilde{C}_\alpha$. 

Therefore, in either case, $H$ is Morse-Kirwan near $p$ with index $0$  if $\lambda_1>0$, and index $2$ if $\lambda_1<0$. The  case where $\frac{\partial (f\circ g)}{\partial q_1}(0)=0$ and $\frac{\partial (f\circ g)}{\partial q_2}(0)\neq0$ is proved analogously. Consequently, $H$ is Morse-Kirwan near $p$ with index 0, 2, or 4 at $p$. \end{proof}

Next, we give a sufficient condition for $H=f\circ\Phi$ to be Morse–Kirwan near elliptic-elliptic critical point $p$, allowing $\Phi(p)$ to be a critical point of $f$.

\begin{theorem}\label{EE MK local min}
Let $\Phi:M\to\mathbb{R}^2$ be an almost-toric system on  $(M,\omega)$, and let
$f:\mathbb{R}^2\to\mathbb{R}$ be a Morse-Kirwan function. Let $p$ be an elliptic-elliptic point of $\Phi$. Suppose:
\begin{itemize}
\item  $\Phi(p)$ is a local minimum of $f$.
\item The restriction $f|_{\Phi(\overline{C})}$ is a Morse-Kirwan function near $\Phi(p)$ for each connected component $C\subset \operatorname{Crit}_1(\Phi)$ such that $p\in \overline{C}$.
\end{itemize}

Then $f\circ\Phi$ is a Morse-Kirwan function near $p$, with index 0 at every critical point sufficiently near $p$.
\end{theorem}
\begin{proof}
Let $H:=f\circ \Phi$ and we may assume $H(p)=f(\Phi(p))=0.$ Since $\Phi(p)$ is a local minimum of $f$, continuity of $\Phi$ implies that  $p$ is a local minimum of $H$. Working in the local coordinates and notation of Theorem \ref{EE MK}, we may assume, after shrinking $U$ if necessary, that $H\geq 0$ on $U$ and $dH = \frac{\partial (f\circ g)}{\partial q_1}(x_1\,dx_1 + \xi_1\,d\xi_1) + \frac{\partial (f\circ g)}{\partial q_2}(x_2\,dx_2 + \xi_2\,d\xi_2)$.

Define $\varphi(q_2):=(f\circ g)(0,q_2)$ and $\psi(q_1):=(f\circ g)(q_1,0)$.  Then, after shrinking $U$, the critical set of $H$ on $U$ is given by 
\[
\mathrm{Crit}(H)\cap U=A\cup B \cup C
\]
where $A:=\Phi^{-1}(\mathrm{Crit}(f))\cap U$, $B:=\{(0,0,x_2,\xi_2)\in U\; | \;  q_2\in \text{Crit}(\varphi)\}$, and $C:=\{(x_1,\xi_1,0,0)\in U\; | \;  q_1\in \text{Crit}(\psi)\}$.

We claim that $A\cup B \cup C\subset H^{-1}(0)$. 
We first show $A\subset H^{-1}(0)$. Since $f$ is Morse-Kirwan and $\Phi(p)\in \text{Crit}(f)$ with $f(\Phi(p))=0$, Proposition \ref{equivalence of locally finitness} shows that there exists a neighborhood $V$ of $\Phi(p)$ such that $\text{Crit}(f)\cap V\subset f^{-1}(0)$. Shrinking $U$ so that $\Phi(U)\subset V$, we have
\[A=\Phi^{-1}(\text{Crit}(f))\cap U \subset \Phi^{-1}(f^{-1}(0))=H^{-1}(0).\]

Now we show $B\subset H^{-1}(0)$. Since the map $q_2\mapsto g(0,q_2)$ parametrizes $\Phi(\overline{C_1})$ near $\Phi(p)$ where $C_1:=\{(0,0,x_2,\xi_2)\in U \; | \; (x_2,\xi_2)\neq(0,0)\}$, the hypothesis that $f|_{\Phi(\overline{C_1})}$ is Morse-Kirwan near $\Phi(p)$ implies that $\varphi$ is Morse-Kirwan near 0. Moreover, 
since $\Phi(p)$ is a local minimum of $f$ and
$\varphi(q_2)=f \big(g(0,q_2)\big)\geq f(\Phi(p))=\varphi(0)$, 0 is a local minimum of $\varphi$, thus $0\in \text{Crit}(\varphi)$. By Proposition \ref{equivalence of locally finitness}, there exists a neighborhood $W$ of 0 in $q_2(U)$ such that  $\text{Crit}(\varphi)\cap W \subset \varphi^{-1}(0)$. After shrinking $U$ so that $q_2(U)\subset W$, every $x\in B$ satisfies
\[H(x)=f(g(0,q_2))=\varphi(q_2)=0.\]
Hence $B\subset H^{-1}(0).$ Analogously, we have $C\subset H^{-1}(0)$. Consequently,  $\text{Crit}(H)\cap U=A\cup B \cup C \subset H^{-1}(0)$.

Since  $p$ is a local minimum of $H$ and $\text{Crit(H)}\cap U \subset H^{-1}(0)$, Lemma \ref{local min and preimage} implies that $H$ is Morse-Kirwan near $p$. 
Moreover, every $x\in \text{Crit}(H)\cap U$ is a local minimum of $H$ since  $\text{Crit}(H)\cap U\subset H^{-1}(0)$ and $H\geq0$ on $U$.  Thus $\text{Hess}_x(H)$ is positive semidefinite for every $x\in \text{Crit}(H)\cap U$, and hence $H$ has index $0$ on $\text{Crit}(H)\cap U$. \end{proof}

\section{Local Analysis near Focus-Focus Singularities}\label{sec:FF}

In this section, we study the composition $H=f\circ \Phi$ near focus-focus singularity $p$. We first consider the case where $\Phi(p)$ is not a critical point of $f$. 

\begin{theorem}\label{FF MK}
Let $\Phi:M\to\mathbb{R}^2$ be a nondegenerate integrable system on $(M,\omega)$, and let $f:\mathbb{R}^2\to\mathbb{R}$ be a smooth function. Let $p$ be a focus-focus critical point of $\Phi$, and suppose  that $\Phi(p)\notin \operatorname{Crit}(f)$. Then $f\circ\Phi$ is a Morse  function  near $p$, with index 2 at $p$.
\end{theorem}

\begin{proof}
By Theorem \ref{local normal form theorem}, there exist local symplectic coordinates $(x_1,\xi_1,x_2,\xi_2)$ centered at $p$  on a sufficiently small neighborhood $U$, and a local diffeomorphism $g:\mathbb{R}^2\rightarrow \mathbb{R}^2$ with $g(0)=\Phi(p)$ such that $\Phi=g\circ q$ where $q=(q_1,q_2)=(x_1\xi_2 - x_2\xi_1  , \; x_1\xi_1 + x_2\xi_2)$. Let $H:=f\circ\Phi=(f\circ g)(q_1,q_2)$. Then the differential of $H$ is
\[
dH
    =\bigg[\frac{\partial (f\circ g)}{\partial q_1}\ \ \frac{\partial (f\circ g)}{\partial q_2}\bigg]\begin{bmatrix}\xi_2 & -x_2 & -\xi_1 & x_1\\
            \xi_1 & x_1 & \xi_2 & x_2\end{bmatrix}.
\]

Since $\Phi(p)\notin \text{Crit}(f)$ and $dg_0$ is invertible, we have $d(f\circ g)_0\neq 0$. Thus, after shrinking $U$ if necessary, $d(f\circ g)\neq 0$ on $q(U)$. It follows that  
\[\text{Crit}(H)\cap U=\{p\}.\]

At $p$  we have $dq_1(p)=dq_2(p)=0$. Therefore, we have 
\[
\text{Hess}_p(H)
=
\sum_{i,j=1}^{2}
\frac{\partial^2(f\circ g)}{\partial q_i\,\partial q_j}(0)\,
dq_i(p)\otimes dq_j(p)
+
\sum_{i=1}^{2}
\frac{\partial(f\circ g)}{\partial q_i}(0)\,\text{Hess}_p(q_i)=\sum_{i=1}^{2}
\frac{\partial(f\circ g)}{\partial q_i}(0)\,\text{Hess}_p(q_i).
\] The eigenvalues of $\text{Hess}_p(H)$ are $\lambda,\lambda,-\lambda,-\lambda,$  where $\lambda=\sqrt{\big(\frac{\partial (f\circ g)}{\partial q_1}(0)\big)^2+\big(\frac{\partial (f\circ g)}{\partial q_2}{(0)}\big)^2}>0. $ Therefore $p$ is a nondegenerate critical point of $H$ of index 2. Since $p$ is the only critical point of $H$ in $U$, $H$ is Morse near $p$.\end{proof}

Next, we give a sufficient condition for $H=f\circ \Phi$ to be Morse--Kirwan near a focus-focus critical point $p$, allowing $\Phi(p)$ to be a critical point of $f$.

\begin{theorem}\label{FF MK local minimum}
Let $\Phi:M\to\mathbb{R}^2$ be a nondegenerate integrable system on  $(M,\omega)$, and let $f:\mathbb{R}^2\to\mathbb{R}$ be a Morse--Kirwan function. Let $p$ be a focus-focus critical point of $\Phi$, and suppose $\Phi(p)$ is a local minimum of $f$. Then $f\circ\Phi$ is a Morse--Kirwan function  near $p$, with index 0 at every critical point sufficiently near $p$.
\end{theorem}

\begin{proof}
By Theorem \ref{local normal form theorem}, $p$ is an isolated critical point of $\Phi$ in a sufficiently small neighborhood $U$. Hence, by Lemma \ref{either a regular or an isolated point}, $H:=f\circ\Phi$ is Morse-Kirwan near $p$. Moreover, after shrinking $U$ if necessary, we have  $\text{Crit}(H)\cap U=\Phi^{-1}(\text{Crit}(f))\cap U.$

Without loss of generality, we may assume  $H(p)=f(\Phi(p))=0.$ Since $\Phi(p)$ is a local minimum of $f$, after shrinking $U$ if necessary, we have $H\geq 0$ on $U$.  Since $f$ is Morse–Kirwan and $\Phi(p)$ is a critical point of $f$ with $f(\Phi(p))=0$, Proposition \ref{equivalence of locally finitness} gives a neighborhood $V$ of $\Phi(p)$ such that $\mathrm{Crit}(f)\cap V \subset f^{-1}(0).$ After further shrinking $U$ so that $\Phi(U)\subset V$, we have  \[\text{Crit(H)}\cap U=\Phi^{-1}(\text{Crit}(f))\cap U  \subset \Phi^{-1}(f^{-1}(0))= H^{-1}(0).\] Thus every $x\in \text{Crit(H)}\cap U$ is a local minimum of $H$. Hence $\text{Hess}_x(H)$ is positive semidefinite, and therefore $H$ has index 0 on $\text{Crit(H)}\cap U$. \end{proof}

\section{Morse-Kirwan compositions}\label{sec:6}

In this section, we combine the local results of Sections \ref{sec:ER}–\ref{sec:FF} to establish sufficient conditions under which \(H=f\circ\Phi\) is a Morse–Kirwan function on \(M\). Then we  give additional conditions that ensure \(H\) is proper, bounded below, and has no critical points of index \(1\).

\begin{theorem}\label{main MK}
Let $\Phi:M\to\mathbb{R}^2$ be an almost-toric system on  $(M,\omega)$, and let
$f:\mathbb{R}^2\to\mathbb{R}$ be a Morse-Kirwan function.  Suppose the following conditions hold:

\begin{itemize}
\item If $p \in \operatorname{Crit}(\Phi)$ and  $\Phi(p)\in \operatorname{Crit}(f)$, then $\Phi(p)$ is a local minimum of $f$.
\item For each connected component $C\subset \operatorname{Crit}_1(\Phi)$ with closure $\overline{C}$,
\begin{itemize}
 \item The restriction $f|_{\Phi(\overline{C})}$ is a Morse-Kirwan function.
    \item If $p\in \overline{C}\setminus C$, then $\Phi(p)$ is either a regular point or a local minimum of  $f|_{\Phi(\overline{C})}$.
   
\end{itemize}

\end{itemize}

Then $f\circ\Phi$ is a Morse-Kirwan function.
\end{theorem}

\begin{proof}
Let \(H:=f\circ\Phi\). Since $dH_p=d f_{\Phi(p)}\circ d\Phi_p$, we have $p\in\operatorname{Crit}(H)$ if and only if $\operatorname{Im}(d\Phi_p)\subset \ker(d f_{\Phi(p)})$.  Thus,
\[
\operatorname{Crit}(H)
=
\Phi^{-1}(\operatorname{Crit}(f))
\sqcup
\left\{
p\in\operatorname{Crit}(\Phi):
\Phi(p)\notin\operatorname{Crit}(f),\
\operatorname{Im}(d\Phi_p)\subset\ker(d f_{\Phi(p)})
\right\}.
\]

Fix \(p\in\operatorname{Crit}(H)\). Suppose first that  \(p\) is a regular point of
\(\Phi\). Then \(\Phi\) is a submersion near \(p\), so  Proposition \ref{submersion and MK}
implies that \(H\) is Morse--Kirwan near \(p\).

Now suppose that  \(p\in\operatorname{Crit}(\Phi)\). Since \(\Phi\) is almost-toric, \(p\) is of elliptic--regular, elliptic--elliptic, or focus-focus. If \(\Phi(p)\notin\operatorname{Crit}(f)\), then
Theorems \ref{ER MK}, \ref{EE MK}, and \ref{FF MK} show that \(H\) is Morse--Kirwan near \(p\) in the three respective cases. If \(\Phi(p)\in\operatorname{Crit}(f)\), then, by hypothesis,
\(\Phi(p)\) is a local minimum of \(f\). Hence Theorems \ref{ER MK local min}, \ref{EE MK local min}, and \ref{FF MK local minimum} show that \(H\) is Morse--Kirwan near \(p\) in the three respective cases. 

Therefore, $H$ is Morse--Kirwan near every critical point of $H$. Since $H$ is automatically Morse--Kirwan near every regular point, Theorem \ref{local to global} imples that $H$ is Morse--Kirwan on \(M\).

\end{proof}

\begin{theorem}\label{main theorem in section 6}
   Let $\Phi:M\to\mathbb{R}^2$ be a proper almost-toric  system on  $(M,\omega)$, and let $f:\mathbb{R}^2\to\mathbb{R}$ be a proper, bounded below Morse-Kirwan function. Suppose the following conditions hold:

\begin{itemize}
\item $f$ has no critical points of index 1 in $\Phi(M)$.
\item If $p \in \operatorname{Crit}(\Phi)$ and  $\Phi(p)\in \operatorname{Crit}(f)$, then $\Phi(p)$ is a local minimum of $f$.
\item For each connected component $C\subset \operatorname{Crit}_1(\Phi)$ with closure $\overline{C}$,
\begin{itemize} 
\item The restriction $f|_{\Phi(\overline{C})}$ is a Morse-Kirwan function.
\item If $p\in \overline{C}\setminus C$, then $\Phi(p)$ is either a regular point or a local minimum of  $f|_{\Phi(\overline{C})}$.
\item If $p\in C$ and $\Phi(p)$ is a local maximum of $f|_{\Phi(C)}$, then the gradient $\nabla f_{\Phi(p)}$ points locally outward from  $\Phi(M)$.   
\end{itemize}

\end{itemize}

Then $f\circ\Phi$ is a proper, bounded below Morse-Kirwan function with no points of index 1.
\end{theorem}
\begin{proof}
Let $H:=f\circ \Phi$. By theorem \ref{main MK}, $H$ is Morse-Kirwan. Since $f$ is bounded below, so is $H$. In addition, since both $f$ and $\Phi$ are proper, $H$ is also proper. It remains to show that $H$ has no critical points of index 1. 

Fix $p\in \text{Crit}(H)$. Suppose first that $p$ is a regular point of $\Phi$. Since $d\Phi_p$ is surjective,  $dH_p=0$ if and only if $df_{\Phi(p)}=0$. Thus $\Phi(p)\in \text{Crit}(f)$, so the Hessian chain rule gives \[\text{Hess}_p(H)(v,w)=\text{Hess}_{\Phi(p)}(f)\big(d\Phi_p(v),d\Phi_p(w) \big)\] for $v,w\in T_pM$. Hence the Hessian of $H$ transverse to $\text{ker}(d\Phi_p)$ agrees with $\text{Hess}_{\Phi(p)}(f)$, while it vanishes on $\text{ker}(d\Phi_p)$. Therefore, $\text{index}_p(H)=\text{index}_{\Phi(p)}(f)$. By hypothesis, $f$ has no  critical points of index 1 in $\Phi(M)$, so $\text{index}_p(H)\neq 1$.

    Now suppose that  $p\in \text{Crit}(\Phi)$. If $\Phi(p)\notin \text{Crit}(f)$,  theorem \ref{ER MK}, \ref{EE MK} and \ref{FF MK} show that $H$ does not have index 1 at $p$. If  $\Phi(p)\in \text{Crit}(f)$, theorem \ref{ER MK local min}, \ref{EE MK local min}, and \ref{FF MK local minimum} show that $H$ does not have index 1 at $p$. Therefore, $H$ is a proper bounded below Morse-Kirwan function with no critical points of index 1.
\end{proof}

\section{Proofs of the main results}\label{sec:main}

In this section, we prove our main results (Theorem \ref{theorem 1}--\ref{theorem 3}) by applying the Morse--Kirwan criteria from  Section \ref{sec:6} to the norm-squared function $f(x,y)=\|(x,y)-(a,b)\|^2$. We first establish two sufficient conditions for the connectedness of a single fiber, which will be used to prove our main theorems.

\begin{proposition}\label{main theorem}
   Let $\Phi:M\to\mathbb{R}^2$ be a proper almost-toric system on $(M,\omega)$. Fix  $(a,b)\in \Phi(M)$, and define $f:\mathbb{R}^2\to\mathbb{R}$ by $f(x,y)=\|(x,y)-(a,b) \|^2$. For each connected component $C\subset Crit_1(\Phi)$ with closure $\overline{C}$, suppose that the following conditions hold:
   \begin{itemize} 
\item The restriction $f|_{\Phi(\overline{C})}$ is a Morse-Kirwan function.
\item If $p\in \overline{C}\setminus C$, then $\Phi(p)$ is either a regular point or a local minimum of  $f|_{\Phi(\overline{C})}$.
\item If $p\in C$ and $\Phi(p)$ is a local maximum of $f|_{\Phi(C)}$, then the gradient $\nabla f_{\Phi(p)}$ points locally outward from  $\Phi(M)$.   
\end{itemize}

Then $\Phi^{-1}(a,b)$ is connected.
\end{proposition}
\begin{proof}
Let $H:=f\circ \Phi=\|\Phi-(a,b)\|^2$. Then $H\geq 0$ and $H^{-1}(0)=\Phi^{-1}(a,b)$. Note that $f(x,y)=\|(x,y)-(a,b) \|^2$ is proper, bounded below, and Morse with  a unique critical point $(a,b)$ of index 0. Moreover, if $p\in\text{Crit}(\Phi)$ and $\Phi(p)\in \text{Crit}(f)$, then $\Phi(p)=(a,b)$, which is the global minimum of $f$. Thus Theorem \ref{main theorem in section 6} implies that $H$ is a proper, bounded below Morse-Kirwan function with no critical points of index 1. 
Then, by Lemma $\ref{main lemma}$, the index 0 critical set of $H$ is connected. Since $H$ is constant on each connected component of its critical set and $H^{-1}(0)$ is a nonempty subset of the index 0 critical set of $H$, it follows that the index 0 critical set of $H$ equals $H^{-1}(0)$. Hence $H^{-1}(0)=\Phi^{-1}(a,b)$ is connected. 

\end{proof}



\begin{corollary}\label{medium main lemma}
   Let $\Phi:M\to\mathbb{R}^2$ be a proper almost-toric system on  $(M,\omega)$. Fix  $(a,b)\in \Phi(M)$, and define $f:\mathbb{R}^2\to\mathbb{R}$ by $f(x,y)=\|(x,y)-(a,b) \|^2$. Suppose that, for each connected component $C\subset \operatorname{Crit}_1(\Phi)$ with closure $\overline{C}$, every critical point of the restriction $f|_{\Phi(\overline{C})}$ is a local minimum. Then $\Phi^{-1}(a,b)$ is connected.

\end{corollary}
\begin{proof}
    Fix a connected component $C\subset \text{Crit}_1(\Phi)$. It suffices to verify the hypotheses of Proposition \ref{main theorem}. Let $\eta:I\rightarrow \Phi(\overline{C})$ be any smooth regular local parametrization. By hypothesis, every critical point of $f\circ \eta$ is a local minimum. Hence, by Proposition \ref{only local min} and Remark \ref{only local min remark},  $f\circ \eta$ is  Morse--Kirwan. Thus $f|_{\Phi(\overline{C})}$ is Morse-Kirwan. Moreover, the remaining hypotheses of Proposition \ref{main theorem} follow from the assumption that every critical point of $f|_{\Phi(\overline{C})}$ is a local minimum.
    Consequently, by Proposition \ref{main theorem}, $\Phi^{-1}(a,b)$ is connected.

\end{proof}


We now apply Corollary \ref{medium main lemma} to prove Theorem \ref{theorem 1}, which gives a sufficient condition for all fibers to be connected when the elliptic-regular critical value branches  lie in affine lines, possibly after a diffeomorphism of the image.

\begin{proof}[Proof of Theorem \ref{theorem 1}]
Let $\widetilde{\Phi}:=g\circ\Phi.$
Since \(\Phi\) is proper and \(g\) is a diffeomorphism onto its closed image $g\big(\Phi(M) \big)$ in $\mathbb{R}^2$,  \(\widetilde{\Phi}\) is proper. In particular, $\operatorname{Crit}_1(\widetilde{\Phi})
=\operatorname{Crit}_1(\Phi).$

    Fix $(a,b)\in \widetilde{\Phi}(M)$ and define $f(x,y):=\|(x,y)-(a,b)\|^2$. Fix a connected component $C\subset \text{Crit}_1(\widetilde{\Phi})$. By hypothesis, $\widetilde{\Phi}(C)$ is contained in an affine line $L\subset\mathbb{R}^2$. Since $\widetilde{\Phi}$ is continuous and $L$ is closed, $\widetilde{\Phi}(\overline{C})\subseteq \overline{\widetilde{\Phi}(C)}\subseteq L$. Hence $\widetilde{\Phi}(\overline{C})$ is also contained in $L$. 
    
    We claim that every critical point of $f|_{\widetilde{\Phi}(\overline{C})}$ is a local minimum. First, let $p\in C$. Since $\widetilde{\Phi}(C)$ is contained in the affine line $L$, near any $p\in C$ we may parametrize $\widetilde{\Phi}(C)$ by $\gamma(t)=\widetilde{\Phi}(p)+tv$, $t\in (-\epsilon,\epsilon)$   for some $v\neq 0.$ If $\widetilde{\Phi}(p)$ is a critical point of $f|_{\widetilde{\Phi}(C)}$, then $(f\circ\gamma)'(0)
=2\langle\widetilde{\Phi}(p)-(a,b),v\rangle=0.$ Thus,  $(f\circ\gamma)(t)-(f\circ\gamma)(0)
=2t\langle \widetilde{\Phi}(p)-(a,b),v \rangle+t^2\|v\|^2 =t^2\|v\|^2>0$ for sufficiently small $t\neq0$. Hence $\widetilde{\Phi}(p)$ is a local minimum of
$f|_{\widetilde{\Phi}(\overline C)}$. 

Now let $p\in\overline C\setminus C$. Since $\widetilde{\Phi}(\overline{C})$ is also contained in the affine line $L$, the same argument using a local parametrization 
 $\gamma(t)=\widetilde{\Phi}(p)+tv$, $t\in[0,\epsilon)$, shows that whenever $\widetilde{\Phi}(p)$ is a critical point of $f|_{\widetilde{\Phi}(\overline C)}$, it is a local minimum. Thus every critical point of $f|_{\widetilde{\Phi}(\overline C)}$ is a local minimum. Therefore, Corollary \ref{medium main lemma} implies that $\widetilde{\Phi}^{-1}(a,b)$ is connected. Since $(a,b)\in\widetilde{\Phi}(M)$ was arbitrary, all fibers of $\widetilde{\Phi}$ are
connected. Finally, since \(g\) is injective,  for every \(q\in\Phi(M)\),
\[
\widetilde{\Phi}^{-1}(g(q))
=(g\circ\Phi)^{-1}(g(q))
=\Phi^{-1}(q).
\] Hence all fibers of $\Phi$ are connected.
\end{proof}

Examples covered by Theorem \ref{theorem 1} include proper almost-toric systems $\Phi$ whose images are a closed half-plane or a closed infinite strip, provided that the image of the closure of each connected component of $\operatorname{Crit}_1(\Phi)$ is contained in an affine line. The image may also contain focus-focus critical values in its interior.

Although these cases are not directly covered by \cite[Theorems 1--3]{pelayo2015fiber}, their fiber connectivity can also be proved using \cite[Proposition 4.5, Theorem 4.7, and Proposition 5.8]{pelayo2015fiber}, after an appropriate diffeomorphism of the image and with the choice  $f(x,y)=x$. In contrast, our approach uses the norm-squared function and applies directly under the affine-line hypothesis.

Moreover, Theorem \ref{theorem 1} applies to proper almost-toric systems with only focus-focus singularities, since $\operatorname{Crit}_1(\Phi)=\emptyset$ and hence its hypothesis is vacuous. Thus, such systems have connected fibers. In this case, the image of $\Phi$ is necessarily $\mathbb R^2$.\\

Now we  prove Theorem \ref{theorem 2} by establishing fiber connectivity over a dense subset of regular values using Proposition \ref{main theorem}, and then extending it to all fibers using Ehresmann's fibration theorem and Lemma \ref{detour lemma}.

\vspace{3mm}


\begin{proof}[Proof of Theorem \ref{theorem 2}]
Since \(\Phi\) is proper and \(g\) is a diffeomorphism onto its closed image $g\big(\Phi(M) \big)$ in $\mathbb{R}^2$,  $\widetilde{\Phi}=g\circ\Phi$ is proper. In particular, $\operatorname{Crit}_1(\widetilde{\Phi})
=\operatorname{Crit}_1(\Phi).$
Let $B:=\widetilde{\Phi}(M)$, let $B_r$ denote the set of regular values of $\widetilde{\Phi}$, and let $\mathcal{F}$ denote the set of focus-focus critical values of $\widetilde{\Phi}$. We prove this theorem by verifying the hypotheses of Lemma \ref{detour lemma}. 

\vspace{1mm}

Since every elliptic-regular critical value of $\widetilde{\Phi}$ lies in $\partial B$, we have $B_r=\interior(B)\setminus \mathcal{F}$. Hence, by Lemma \ref{locally connectedness of B_r}, $B_r\cap D$ is connected for every sufficiently small disc $D$ centered at a point of $B$.  Also, $\widetilde{\Phi}$ is proper and $\text{Crit}(\widetilde{\Phi})$ has empty interior. Therefore  it remains to show that every regular fibers of $\widetilde{\Phi}$ is connected. To prove this, we first construct a dense subset $\mathcal{R}$ of $B_r$ by excluding three measure zero sets and apply Proposition \ref{main theorem} to show that every fiber over $\mathcal{R}$ is connected. Then we extend this connectedness to all regular fibers using Ehresmann’s fibration theorem. 

\vspace{1mm}

Choose countably many regular local parametrizations $\gamma_j:I_j \rightarrow\mathbb{R}^2$ covering the elliptic-regular critical values of $\widetilde{\Phi}$, and  let $n_j(t)$ be the inward unit normal vector along $\gamma_j$. Then the curvature of $\gamma_j$ with respect to $n_j$ is  $\kappa_j(t):=\frac{\langle \gamma_j''(t), \ n_j(t) \rangle}{\|\gamma_j'(t)\|^2}$. Since $\gamma_j$ lies in the boundary of the convex set $B$, we have $\kappa_j(t)\geq0$. At every $t$ with $\kappa_j(t)>0$, the osculating circle of $\gamma_j$ has radius $\frac{1}{\kappa_j(t)}$ and center 
\[O_j(t):=\gamma_j(t)+\frac{1}{\kappa_j(t)}n_j(t).\] Thus $O_j$ defines a  smooth map $\{t\in I_j \ 
 | \ \kappa_j(t)>0\} \rightarrow \mathbb{R}^2 $ assigning to each point the center of its osculating circle. Fix  $c\in \interior{(B)}$ and define $f_c(x,y)=\|(x,y)-c\|^2$.  At a critical point $t_0$ of $f_c\circ \gamma_j=\|\gamma_j-c\|^2$, we have $0=(f_c\circ \gamma_j)'(t_0)=-2 \langle c-\gamma_j(t_0), \gamma_j'(t_0)\rangle$. Therefore $c-\gamma_j(t_0)$ is perpendicular to $\gamma_j'(t_0)$, hence $c-\gamma(t_0)=r\,n(t_0)$ where $r=\|c-\gamma(t_0) \|>0$. Thus, \[(f_c \circ \gamma_j)''(t_0)=2\|\gamma_j'(t_0) \|^2-2r\langle n_j(t_0),\gamma_j''(t_0)\rangle=2\|\gamma_j'(t_0)\|^2\big(1-r\kappa_j(t_0) \big).\] It follows that $t_0$ is a degenerate critical point of $f_c\circ \gamma_j$ if and only if $c=O_j(t_0)$ and $\kappa_j(t_0)\neq0$. If $\kappa_j(t_0)=0$, we have $(f_c \circ \gamma_j)''(t_0)>0$, so degeneracy is impossible in this case. Since  each $O_j$ is a smooth map from an open subset of $\mathbb{R}$ to $\mathbb{R}^2$, Sard's theorem implies that its image has  Lebesgue measure zero. Therefore the countable union $\mathcal{O}:=\bigcup_jO_j\big( \{t\in I_j \ | \ \kappa_j(t)>0 \} \big) $ has measure zero.

\vspace{1mm}

    Now, at each elliptic-elliptic critical point $p$ of $\widetilde{\Phi}$, let $v_{p,1}, v_{p,2}\neq0 $ be tangent vectors at $\widetilde{\Phi}(p)$ to the two local elliptic-regular critical value branches meeting at $\widetilde{\Phi}(p)$. The ellipitc-elliptic critical value $\widetilde{\Phi}(p)$ is critical for the restriction of $f_c$ to the $i$-th branch, $i=1,2$, exactly when $\langle \widetilde{\Phi}(p)-c, v_{p,i} \rangle=0$. Thus, the corresponding centers $c$ lie on the affine line 
    \[L_{p,i}:=\{c\in \mathbb{R}^2 \ | \ \langle \widetilde{\Phi}(p)-c, v_{p,i} \rangle =0 \},\] which has measure zero in $\mathbb{R}^2$. Since there are only countably many elliptic-elliptic critical points, $\mathcal{L}:=\bigcup_{p,i}L_{p,i}$ also has measure zero.  
    
\vspace{1mm}
    
    Consequently, \[\mathcal{R}:=\interior(B)\setminus(\mathcal{O}\cup \mathcal{L}\cup \mathcal{F})\] is dense in $\interior(B)$, and every  $q\in \mathcal{R}$ is a regular value of $\widetilde{\Phi}$. Fix $q\in \mathcal{R}$.  For each connected component $C\subset \text{Crit}_1(\widetilde{\Phi})$, $q\notin\mathcal{O}$ implies that every critical point of $f_q|_{\widetilde{\Phi}(C)}$ is nondegenerate. Moreover, for every $p\in \overline{C}\setminus C$, $q\notin \mathcal{L}$ implies that $\widetilde{\Phi}(p)$ is a regular point of $f_q|_{\widetilde{\Phi}(\overline{C})}$. Hence, every critical point of $f_q|_{\widetilde{\Phi}(\overline{C})}$ lies in $\widetilde{\Phi}(C)$ and is nondegenerate. Therefore $f_q|_{\widetilde{\Phi}(\overline{C})}$ is Morse, and hence Morse-Kirwan.  Moreover,  if $p\in C$ and $\widetilde{\Phi}(p)$ is a local maximum of $f_q|_{\widetilde{\Phi}(\overline{C})}$,    the gradient $\nabla f_q(\widetilde{\Phi}(p))=2(\widetilde{\Phi}(p)-q)$ points outward from $B$, since $q\in \interior(B)$ and $B$ is convex. Thus Proposition \ref{main theorem} implies that $\widetilde{\Phi}^{-1}(q)$ is connected. Since $q\in \mathcal{R}$ was arbitrary, $\widetilde{\Phi}^{-1}(q)$ is connected for every $q\in \mathcal{R}$. 

    Let $q_0\in B_r$ be an arbitrary regular value of $\widetilde{\Phi}$. Since  $\widetilde{\Phi}$ is proper, Ehresmann's fibration theorem implies that there exists a neighborhood $U_{q_0}\subset B_r$ of $q_0$ such that  $\widetilde{\Phi}|_{\widetilde{\Phi}^{-1}(U_{q_0})}:\widetilde{\Phi}^{-1}(U_{q_0})\rightarrow U_{q_0}$ is a locally trivial fibration. Thus all fibers over $U_{q_0}$ are diffeomorphic.  Since  $\mathcal{R}$ is dense in $B_r$, there exists $q\in U_{q_0}\cap \mathcal{R}$, and hence $\widetilde{\Phi}^{-1}(q_0)$ is diffeomorphic to  $\widetilde{\Phi}^{-1}(q)$. Since $\widetilde{\Phi}^{-1}(q)$ is connected, $\widetilde{\Phi}^{-1}(q_0)$ is also connected. Therefore, every regular fiber of $\widetilde{\Phi}$ is connected.
    
    Consequently, Lemma \ref{detour lemma} implies that all fibers of $\widetilde{\Phi}$ are connected.  Since \(g\) is injective, for every \(q\in\Phi(M)\),  \[\widetilde{\Phi}^{-1}(g(q))
=(g\circ\Phi)^{-1}(g(q))
=\Phi^{-1}(q).\] Hence all fibers of $\Phi$ are connected.  \end{proof}

Theorem \ref{theorem 2} generalizes \cite[Theorem 3]{pelayo2015fiber} in two ways. First, it recovers the fiber connectivity results for all four cases considered in \cite[Theorem 3]{pelayo2015fiber}. Second, it extends these results to noncompact symplectic manifolds with closed convex moment map images, provided that the almost-toric system is proper.\\

Now we prove Theorem \ref{theorem 3}. Unlike in the proof of Theorem \ref{theorem 2}, we can apply Corollary \ref{medium main lemma} to show the connectedness of a single regular fiber by choosing the center of the norm-squared function sufficiently close to the origin. Then we extend this connectedness to  all fibers using Ehresmann's fibration theorem and Lemma \ref{detour lemma}.

\vspace{3mm}

\begin{proof}[Proof of Theorem \ref{theorem 3}]
Since \(\Phi\) is proper and \(g\) is a diffeomorphism onto its closed image $g\big(\Phi(M) \big)$ in $\mathbb{R}^2$,  $\widetilde{\Phi}=g\circ\Phi$ is proper.
We first show that one regular fiber of $\widetilde{\Phi}$ is connected using Corollary \ref{medium main lemma}. Then we extend this connectedness to all regular  fibers using Ehresmann's fibration theorem, and finally apply Lemma \ref{detour lemma} to conclude that all fibers are connected.

Each branch $\Gamma_i$  has the parametrization $\gamma_i(t)=R_i\big(\alpha\sqrt{1+t^2},\beta t \big), \ t\in \mathbb{R}$. Choose $0<\varepsilon<\alpha$ such that $D_{\varepsilon}(0)\subset \text{Int}(\widetilde{\Phi}(M))$. Fix $c\in D_{\varepsilon}(0)$, and define $f_c(x,y)=\|(x,y)-c\|^2$. For each $i$, let $(a_i,b_i)=R_i^{-1}(c)$. Since rotation preserves the Euclidean norm, we have $|a_i|\leq \sqrt{a_i^2+b_i^2}=\|c\|<\alpha$.

Since rotations preserve distances, 
\begin{align*}
    h_{i}(t):=(f_c\circ\gamma_i)(t)&=\|\gamma_i(t)-c\|^2\\
    &=\|R_i^{-1}\big( \gamma_i(t)-c \big)\|^2\\
    &=\| \big(\alpha\sqrt{1+t^2},\beta t \big) - R_i^{-1}(c)    \|^2\\
    &=\big(\alpha \sqrt{1+t^2} -a_i \big)^2 +(\beta t -b_i)^2.
\end{align*}
Differentiating twice  gives 
\[h_i''(t)=2(\alpha^2 +\beta^2)-\frac{2\alpha a_i}{(1+t^2)^{\frac{3}{2}}}\geq 2(\alpha^2+\beta^2-\alpha\|c\|)>2\beta^2>0.\] 
Thus every critical point of $f_c|_{\Gamma_i}$ is a local minimum. Therefore, Corollary \ref{medium main lemma} implies that $\widetilde{\Phi}^{-1}(c)$ is connected for every $c\in D_\varepsilon(0)$. 

Let $\mathcal{F}$ denote the set of focus-focus critical values  of $\widetilde{\Phi}$, and let $B_r$ denote its set of regular values. Since $D_\varepsilon(0)\subset \text{Int}(\widetilde{\Phi}(M))$, the only possible critical values in $D_\varepsilon(0)$ are focus-focus values. As $\mathcal{F}$ is locally finite, we may choose  $c_*\in D_\varepsilon(0)\cap B_r$. Then $\widetilde{\Phi}^{-1}(c_*)$ is a connected regular fiber.

By the boundary assumption, we have $B_r=\text{Int}(\widetilde{\Phi}(M))\setminus\mathcal{F}$. Since $\widetilde{\Phi}(M)$ is connected and its boundary consists of the $n$ disjoint hyperbolic branches $\Gamma_i$, $\text{Int}\big(\widetilde{\Phi}(M)\big)$ is connected. Moreover, since $\mathcal{F}$ is locally finite, $B_r$ is connected. Since $\widetilde{\Phi}$ is proper, Ehresmann's fibration theorem implies that  $\widetilde{\Phi}|_{\widetilde{\Phi}^{-1}(B_r)}:\widetilde{\Phi}^{-1}(B_r)\rightarrow B_r$ is a locally trivial fibration. Since $B_r$ is connected, all regular fibers are diffeomorphic. In particular, since  $\widetilde{\Phi}^{-1}(c_*)$ is connected, every regular fiber of $\widetilde{\Phi}$ is connected. 

Now we verify that $B_r\cap D$ is connected for every critical value $q$ of $\widetilde{\Phi}$ and for every sufficiently small disk $D$ centered at $q$. If $q$ is a focus-focus critical value, then  $q\in\text{Int}(\widetilde{\Phi}(M))$ and $q$ is isolated. Thus, for every sufficiently small disk $D$ centered at $q$, we have $D\subset \text{Int}(\widetilde{\Phi}(M))$ and $D\cap \mathcal{F}=\{q\}$. Consequently, $B_r\cap D=D\setminus \{q\}$ is connected.

If $q$ is an elliptic-regular critical value, then $q\in \Gamma_i$ for some $i$. Since the boundary branches are pairwise disjoint smooth embedded curves, for every sufficiently small disk $D$ centered at $q$, we have $D\cap \partial\widetilde{\Phi}(M)=D\cap \Gamma_i$, and $D\cap \Gamma_i$ is a single embedded arc dividing $D$ into two connected open components. By the elliptic-regular local normal form, $\text{Int}(\widetilde{\Phi}(M))\cap D$ is exactly one of these two components, therefore $\text{Int}(\widetilde{\Phi}(M))\cap D$ is connected. Shrinking $D$ further so that $D\cap \mathcal{F}=\emptyset$, we obtain $B_r\cap D=\big( \text{Int}(\widetilde{\Phi}(M))\setminus \mathcal{ F}\big)\cap D= \text{Int}(\widetilde{\Phi}(M))\cap D$, which is connected. Thus the local connectedness hypothesis of Lemma \ref{detour lemma} holds.

Consequently, since $\widetilde{\Phi}$ is proper, all regular fibers are connected, and $\text{Crit}(\widetilde{\Phi})$ has empty interior, Lemma \ref{detour lemma} implies that all fibers of $\widetilde{\Phi}$ are connected. Since \(g\) is injective, for every \(q\in\Phi(M)\),
\[
\Phi^{-1}(q)
=
\widetilde{\Phi}^{-1}(g(q)).
\]
Therefore all fibers of \(\Phi\) are connected.\end{proof}

Theorem \ref{theorem 3} provides a further family of noncompact, nonconvex moment map images for which fiber connectivity holds. Unlike the images considered in Theorems \ref{theorem 1} and \ref{theorem 2}, these images are neither convex nor bounded by affine lines.





\bibliographystyle{plain}
\bibliography{ref.bib}

\begin{thebibliography}{10}

\bibitem{atiyah1982convexity}
Michael~Francis Atiyah.
\newblock Convexity and commuting hamiltonians.
\newblock {\em Bulletin of the London Mathematical Society}, 14(1):1--15, 1982.

\bibitem{bolsinov2004integrable}
Aleksei~Viktorovich Bolsinov and Anatolij~Timofeevi{\v{c}} Fomenko.
\newblock {\em Integrable Hamiltonian systems: geometry, topology, classification}.
\newblock CRC press, 2004.

\bibitem{Eli84}
Håkan Eliasson.
\newblock {\em Hamiltonian systems with Poisson commuting integrals}.
\newblock PhD thesis, Stockholm University, 1984.

\bibitem{eliasson1990normal}
LH~Eliasson.
\newblock Normal forms for hamiltonian systems with poisson commuting integrals—elliptic case.
\newblock {\em Commentarii Mathematici Helvetici}, 65(1):4--35, 1990.

\bibitem{GuSt1982}
V.~Guillemin and S.~Sternberg.
\newblock Convexity properties of the moment mapping.
\newblock {\em Invent. Math.}, 67:491--513, 1982.

\bibitem{holm2016morse}
Tara Holm and Yael Karshon.
\newblock The morse--bott--kirwan condition is local.
\newblock {\em Research in the Mathematical Sciences}, 3(1):25, 2016.

\bibitem{kirwan1985cohomology}
FC~Kirwan.
\newblock Cohomology of quotients in algebraic and symplectic geometry.
\newblock {\em Math. Notes}, 31, 1985.

\bibitem{lerman1996non}
Eugene Lerman, Eckhard Meinrenken, Sue Tolman, and Chris Woodward.
\newblock Non-abelian convexity by symplectic cuts.
\newblock {\em arXiv preprint dg-ga/9603015}, 1996.

\bibitem{lerman1997hamiltonian}
Eugene Lerman and Susan Tolman.
\newblock Hamiltonian torus actions on symplectic orbifolds and toric varieties.
\newblock {\em Transactions of the American mathematical society}, 349(10):4201--4230, 1997.

\bibitem{leung2010almost}
Naichung~Conan Leung and Margaret Symington.
\newblock Almost toric symplectic four-manifolds.
\newblock 2010.

\bibitem{miranda2006note}
Eva Miranda and Nguyen~Tien Zung.
\newblock A note on equivariant normal forms of poisson structures.
\newblock {\em Mathematical Research Letters}, 13(6):1001--1012, 2006.

\bibitem{palmer2025integrable}
Joseph Palmer.
\newblock Integrable systems with symmetries: toric, semitoric, and beyond.
\newblock {\em arXiv preprint arXiv:2510.05337}, 2025.

\bibitem{pelayo2015fiber}
{\'A}lvaro Pelayo, Tudor~S Ratiu, et~al.
\newblock Fiber connectivity and bifurcation diagrams of almost toric integrable systems.
\newblock {\em Journal of Symplectic Geometry}, 13(2):343--386, 2015.

\bibitem{pelayo2011symplectic}
{\'A}lvaro Pelayo and San V{\~u}~Ng{\d{o}}c.
\newblock Symplectic theory of completely integrable hamiltonian systems.
\newblock {\em Bulletin of the American Mathematical Society}, 48(3):409--455, 2011.

\bibitem{sepe2025connectedness}
Daniele Sepe and Susan Tolman.
\newblock Connectedness of fibers beyond semitoric systems ii: ephemeral singular points.
\newblock {\em arXiv preprint arXiv:2510.16976}, 2025.

\bibitem{sepe2026connectedness}
Daniele Sepe and Susan Tolman.
\newblock Connectedness of fibers beyond semitoric systems i: the non-degenerate case.
\newblock {\em International Mathematics Research Notices}, 2026(12):rnag129, 2026.

\bibitem{sepe2018integrable}
Daniele Sepe and San V{\~u}~Ng{\d{o}}c.
\newblock Integrable systems, symmetries, and quantization.
\newblock {\em Letters in Mathematical Physics}, 108(3):499--571, 2018.

\bibitem{VN2007}
S.~V\~{u}~Ng\d{o}c.
\newblock Moment polytopes for symplectic manifolds with monodromy.
\newblock {\em Adv. Math.}, 208(2):909--934, 2007.

\bibitem{vu2013smooth}
San V{\~u}~Ng{\d{o}}c and Christophe Wacheux.
\newblock Smooth normal forms for integrable hamiltonian systems near a focus--focus singularity.
\newblock {\em Acta Mathematica Vietnamica}, 38(1):107--122, 2013.

\bibitem{wacheux2013systemes}
C~Wacheux.
\newblock {\em Syst{\`e}mes int{\'e}grables semi-toriques et polytopes moment}.
\newblock PhD thesis, PhD thesis, Universit{\'e} de Rennes 1, 2013.

\bibitem{weinsteinr}
A~Weinstein.
\newblock R. abraham and je marsden. foundations of mechanics. addison-wesley, 1978.(with the assistance of tudor ratiu and richard cushman). 1, 9.

\bibitem{williamson1936algebraic}
John Williamson.
\newblock On the algebraic problem concerning the normal forms of linear dynamical systems.
\newblock {\em American journal of mathematics}, 58(1):141--163, 1936.

\end{thebibliography}

\end{document}